\documentclass[11pt,oneside,reqno]{amsart}
\usepackage{eurosym}
\usepackage{amsmath}
\usepackage{amsfonts}
\usepackage{amsthm}
\usepackage{color}
\usepackage{graphicx}
\usepackage{hyperref}
\usepackage[notref,color]{showkeys}
\usepackage{anysize}
\usepackage{cleveref}

\marginsize{1.5cm}{1.5cm}{1.5cm}{1.5cm}

\DeclareMathOperator{\si}{sgn}
\newtheorem{theorem}{Theorem}
\newtheorem{remark}{Remark}

\newtheorem{lemma}{Lemma}

\begin{document}
\title[Control of surface gravity waves]{Controllability of surface gravity
waves and the sloshing problem}
\author[M. A. Fontelos]{M. A. Fontelos}
\address[M. A. Fontelos]{(Corresponding Author) ICMAT-CSIC, C\slash
Nicol\'as Cabrera, no 13-15 Campus de Cantoblanco, UAM, 28049 Madrid, Spain.}
\email{marco.fontelos@icmat.es}
\author[J. L\'opez-R\'ios]{J. L\'opez-R\'ios}
\address[J. L\'opez-R\'ios]{Universidad Industrial de Santander, Escuela de Matemáticas, A.A. 678, Bucaramanga, Colombia}
\email{jclopezr@uis.edu.co}
\date{\today}
\thanks{M. A. Fontelos was partially supported by Ministerio de Econom\'ia,
Industria y Competitividad, Gobierno de Espa\~na (Grant No.
MTM2017--89423--P)}

\begin{abstract}
We study the problem of controlling the free surface for a two dimensional
solid container in the context of the gravity waves and the sloshing
problem. By using conformal maps and the Dirichlet--Neumann operator, the
problem is formulated as a second order evolutionary equation on the free
surface involving a self-adjoint operator. We present then the appropriate
Sobolev spaces where having solutions for the system and study the exact
controllability through an observability inequality for the adjoint problem.
\end{abstract}

\maketitle

%\begin{keywords}
%	Interior controllability, the sloshing problem, gravity waves, Hilbert transform
%\end{keywords}

%\begin{AMS}
%	76B15, 76B10, 93C95
%\end{AMS}

\section{Introduction}

Finding and controlling the frequencies for the oscillations of the free
boundary in the context of gravity waves and studying related spectral
problems is a classical, well known problem in the literature \cite%
{coron2007control,fontelos2020gravity,godoy2008modeling,zuazua2007controllability}%
. Roughly speaking, controlling a system consist not only in testing that
its behavior is satisfactory, but also in putting things in order to
guarantee that it behaves as desired. In mathematical terms, controlling the
\textit{state} $y$, ruled by the \textit{state equation}
\begin{equation*}
A(y)=f(v),
\end{equation*}
where $v$ is the \textit{control}, consists in finding $v\in U_{ad}$,
\textit{the set of admissible controls}, such that the solution to the
equation gets close to a desired prescribed state, $\bar{y}$.

In this paper, we study controllability of a Partial Differential Equation
(PDE), in the context of controlling the oscillations of a liquid free
surface in a two-dimensional bounded container and the so-called sloshing
problem. We formulate the geometrical problem in terms of an
integrodifferential equation by using the Hilbert transform, then we
establish the appropriate Sobolev spaces to study existence of solutions for
the eigenvalue problem and, finally, we set up an observability inequality
for the homogeneous adjoint problem. The sloshing problem adduces to an
important difference with respect to the classical water-waves formulation:
the presence of vertical walls and the contact with the free surface.
Inspired in the developments for the classical wave equation, we introduce
analytical tools to prove that it is possible to control the oscillations of
the free surface, by injecting fluid on the rigid side walls.

The main strengths of our method lie in the use of the Hilbert transform to
formulate the problem as an evolutionary equation involving a self-adjoint
operator. This is known as the boundary integral method and has proved to be
very fruitful in the study of water waves problems (see \cite%
{fontelos2010singularities} and references therein, for instance). Moreover,
the use of Tchebyshev polynomials provides an explicit orthogonal basis
which allows to study, analytically, the associated eigenvalue problem.
Then, an observability inequality arises naturally.

Generally speaking, the water-waves problem for an ideal liquid consists of
describing the motion of a layer of incompressible, inviscid fluid,
delimited below by a solid bottom, and above by a free surface under the
influence of gravity. In mathematical terms, if $\boldsymbol{u}(x,y)$ is the
fluid velocity and $\varphi$ is the velocity potential such that $%
\boldsymbol{u}=\nabla\varphi$, by the conservation laws \cite%
{lannes2013water}
\begin{equation}  \label{ww}
\left\{ \begin{aligned} &\Delta\varphi=0, && \Omega_t, \\
&\eta_t+\eta_x\varphi_x=\varphi_y, && y=\eta, \\
&\varphi_t+\frac{1}{2}(\varphi_x^2+\varphi_y^2)+g\eta=0, && y=\eta, \\
&\partial_n\varphi=0, && y=b, \end{aligned} \right.
\end{equation}
where $\eta$ and $b$ are the free boundary and bottom parametrization,
respectively, $g$ is the acceleration due to gravity and $\Omega_t=\{(x,y)\in%
\mathbb{R}^2:b(x)<y<\eta(t,x)\}$. We put ourselves in the general situation
described in (\ref{ww}). However, we consider the case of a bounded domain
and the sloshing problem of describing the contact line between the free
surface and the solid walls \cite{fontelos2020gravity}. Two
conditions are customary, the pinned--end boundary condition
where the contact line is always pinned to the solid surface, as considered
in \cite{benjamin1979gravity,graham1983new}, and the free--end condition where the contact angle between the fluid--air
interface and the side walls is fixed and the contact line is allowed to
move \cite{nicolas2005effects}. We will consider both and will
impose conditions on the Cauchy problem to deal with it, accordingly.

Controlling the surface by different methods is of practical interest in
oceanography, controllability and inverse problems theory. We
mention, for example, the work by Reid and Russell \cite{reid1985boundary}
where the authors dealt with the linear conservation laws and the
null-controllability in infinite time of the free surface, by a source
control, in a two dimensional domain with flat side walls. Also, the work by
Reid, \cite{reid1995control}, where the capillary version and the control in
finite time is considered. Concerning nonlinear water-waves, there is the
recent work by Alazard \cite{alazard2017stabilization}, for a two
dimensional rectangular domain, where the stabilization through an external
pressure acting on a small part of the free surface is studied. Also \cite%
{alazard2018boundary}, where the author studied the boundary observability
problem in a three dimensional rectangular domain; namely, an estimate for
the energy of the system in terms of the surface velocity at the contact
line with a vertical wall. Finally, in \cite{alazard2018control}, Alazard et
al. addressed the local exact controllability of the two dimensional full
water-waves system, by controlling a localized portion of the free surface,
through the external pressure. On the literature concerning the generation
of waves by wave-makers, controllability and stability properties in the
water-waves context, we refer to \cite%
{mottelet2000controllability,su2020stabilizability,su2021strong}. From the
optimal control point of view, we mention \cite{nersisyan2014generation},
where the authors designed the `best' moving solid bottom generating a
prescribed wave under the context of a BBM-type equation. We mention also the possibility of studying the inverse problem of detecting the source where jets originate, denoted as $J$, by measuring the free surface as in \cite{lecaros2020stability}.

When we talk about the controllability by fluid injection, for instance, we
mean the condition
\begin{equation*}
\frac{\partial \varphi }{\partial n}=J,
\end{equation*}%
for a given function $J=J(t,x)$, with $n$ being the outward normal vector.
This boundary condition lead to think of a boundary control of the gravity
waves problem; nevertheless, we will restate the problem as an
integrodifferential equation on the free surface and the boundary condition
becomes a source term. Then, we may use the classical approach of interior
controllability by means of the adjoint problem and the observability
inequality \cite{micu2004introduction,zuazua1990introduction}.

By addressing this problem, we give an answer to a practical
question raised in \cite{brimacombe1985toward} and numerically studied in
\cite{godoy2008modeling}; namely, the problem of controlling undesirable
splashing appearing in a cooper converter when air is injected into the
molten matte. In \cite{godoy2008modeling}, the authors studied the problem
by using triangular finite elements to mesh a half-ball bounded domain, on a
damped linear gravity waves model. Our approach allows to
consider any general simply connected two dimensional domain, through a
conformal mapping into the lower half-plane. Moreover, if $f$ represents
such a conformal mapping, the geometry is characterized explicitly by the
term $1/|f^{\prime }|$ appearing as a factor in the evolution problem (see (\ref{cp}) below). On this matter, see \cite{fontelos2020gravity}, where
oscillations are numerically computed for bottoms with rectangular even distributions.

Following the methods introduced in \cite{fontelos2020gravity}, after
linearizing (\ref{ww}), the problem is restated through a conformal map into
the lower half--plane. We link the normal derivative to the specific
conformal map and rewrite the problem as a second order evolution equation
on the interval $[-1,1]$, where the Hilbert transform is involved. In \cite%
{fontelos2020gravity}, we used this approach to propose an efficient
computational method to find the sloshing frequencies on general 2d domains.
We mention \cite{kim2015capillary} where the capillary version of the
problem is developed, under the context of the oscillations in a nozzle of
an inkjet printer. In these works, two possibilities for the contact line
were considered: the `pinned-end edge condition', where the contact line is
always pinned to the solid surface, and the `free-end condition' where the
contact angle between the fluid--air interface and the side walls is fixed
and the contact line is allowed to move, with contact angle $\pi/2$.

With the present work, we have completed the numerical analysis started in
\cite{fontelos2020gravity}. Namely, we establish a Sobolev frame where the
Cauchy problem is well-posed and study the eigenvalue associated problem. We
also prove an observability inequality for the adjoint problem and explore
the possibility of finding explicit controls taking the oscillations at the
free surface to zero.

The rest of the paper is organized as follows: In section \ref{S2} we formulate the general equations to be considered, the linearization approach and corresponding formulation on the half-plane by the conformal mapping. In section \ref{S3} we use the Hilbert transform to state a second order evolutionary
PDE modeling the dynamics of the fluid interface on the bounded domain $%
[-1,1]$. In section \ref{S4} we make use of the Tchebyshev polynomials to study the stationary adjoint problem in suitable Sobolev spaces. In section \ref{S5} we
prove an observability inequality for the adjoint problem. Finally, in section \ref{S6} we establish the controllability of the problem and explore the
possibility of determining possible control functions, explicitly.

\section{Formulation}

\label{S2}

Let us consider a two dimensional container, filled with water, bounded from
above by a free surface. In this context, the motion is governed by the
incompressible Euler equations with zero surface tension:
\begin{align*}
\nabla\cdot\mathbf{u}&=0, \\
\rho\left[\frac{\partial\mathbf{u}}{\partial t}+(\mathbf{u}\cdot\nabla)
\mathbf{u}\right]&=-\nabla p-ge_2,
\end{align*}
with $-ge_2$ being the constant acceleration of gravity, $g>0$, and $e_2$
the unit upward vector in the vertical direction, $\rho$ is the (constant)
density of the fluid and $p$ is the pressure inside. We use the classical
notation $(x,y)\in\mathbb{R}^2$ and $z=x+iy$ for complex numbers.

By considering the potential function $\varphi$ of $\mathbf{u}$, so that $%
\mathbf{u}=\nabla\varphi$:
\begin{align}  \label{pf2}
\Delta\varphi&=0, \\
\frac{\partial\varphi}{\partial t}+\frac{1}{2}|\nabla\varphi|^2+\frac{1}{%
\rho }p+gy&=\text{const}.
\end{align}

We complement system above with the nonzero Neumann condition at the solid
walls, $\mathbf{u}\cdot n=J(t,x)$, and a kinematic condition on the free
boundary. In terms of $\varphi $:
\begin{align}
\frac{\partial \varphi }{\partial n}& =J(t,x),\quad \text{at the solid
boundary }  \label{pf3} \\
\eta _{t}& =\frac{\partial \varphi }{\partial n},\quad |x|<1.  \label{pf4}
\end{align}

\subsection{Linearized equations}

\label{ss1}

Next, we are going to linearize Euler's system around the zero state. This will allow us to obtain a single, explicit, evolutionary conservation law on the free surface modeling the fluid dynamics of the free surface side bounded by vertical solid walls. This equation must be complemented with boundary conditions at the contact line, as explained in section \ref{S3}.

Even if we consider the reference domain $\Omega =\{(x,y)\in \mathbb{D}%
:y^{\prime }<0\}$, other domains can be considered though. Let
\begin{equation*}
\eta (t,x)=\epsilon \zeta (t,x),
\end{equation*}%
and%
\begin{equation*}
J(t,x)=\epsilon j(t,x).
\end{equation*}

Then if
\begin{equation*}
\varphi =const.+\epsilon \phi ,
\end{equation*}%
conditions (\ref{pf2})-(\ref{pf3}) in terms of $\phi $, at the first order
for $\epsilon <<1$, become (after re-scaling to make $g=1$)
\begin{align}
\Delta \phi & =0,  \label{le1} \\
\phi _{t}+\zeta & =0,\quad \text{at }y=0,\ |x|\leq 1,  \label{le2} \\
\zeta _{t}& =\frac{\partial \phi }{\partial n},\quad \text{at }y=0,\ |x|\leq
1,  \label{le3} \\
\frac{\partial \phi }{\partial n}& =j(t,x),\quad \text{on the solid walls.}
\label{le4}
\end{align}

System above is complemented with initial conditions $\phi_0$, $%
\zeta_0$ such that the following mass conservation is satisfied:
\begin{equation}  \label{mass}
\int_{-1}^1\zeta_0(x)dx=0.
\end{equation}

Therefore, from (\ref{le2})-(\ref{le3}), on the free boundary and for $%
|x|\leq 1$,
\begin{equation}
\phi _{tt}+\frac{\partial \phi }{\partial n}=0.  \label{le5}
\end{equation}

\begin{figure}[tbp]
\begin{center}
\includegraphics[scale=0.841]{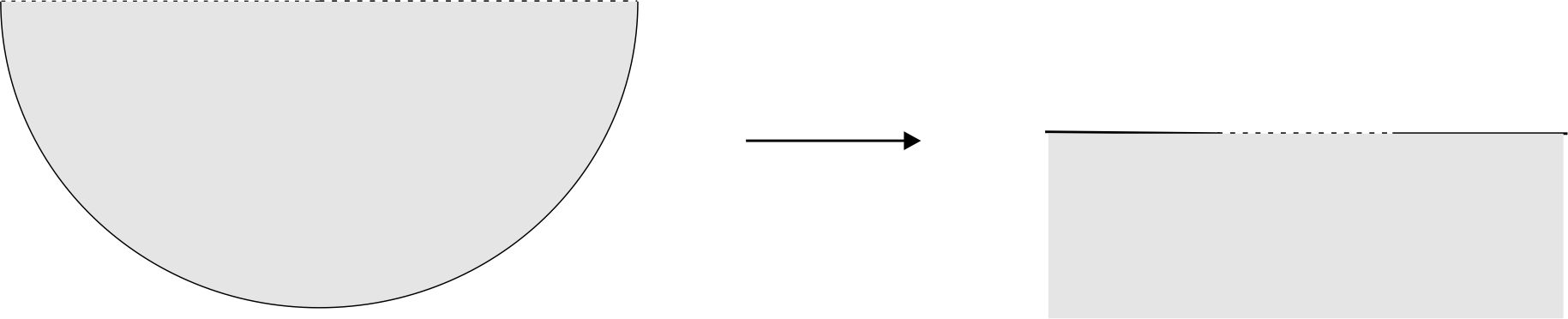} \put(-225,79){$1$} \put(-378,79){$-1$}
\put(-310,40){$\Delta_z\varphi=0$} \put(-350,-6){$\partial_n\varphi=J$}
\put(-220,25){$w=f(z)=\frac{1}{\frac{1}{2}(z+\frac{1}{z})}$} \put(-80,15){$%
\Delta_{w}\widetilde{\varphi}=0$} \put(-130,56){$\partial_{\widetilde{n}}%
\widetilde{\varphi}=\widetilde{J}/|f'|$}
%\put(-40,60){$\partial_{\widetilde{n}}\wph=\widetilde{j}/|f'|$}
\put(-40,35){$1$} \put(-90,35){$-1$}
\end{center}
\caption{Geometry of the problem.}
\label{F1}
\end{figure}
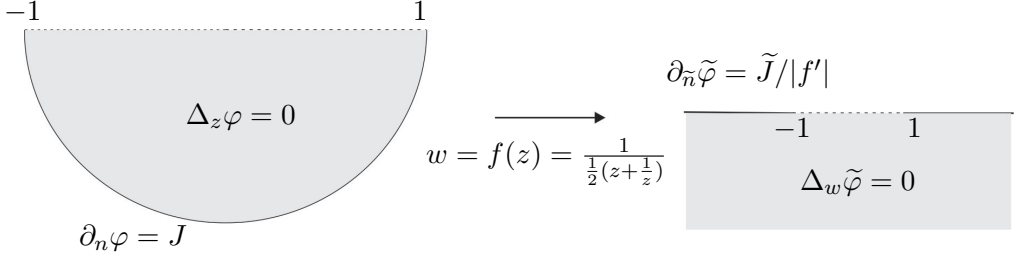

We remark that the motion of the fluid interface may also be affected by
external forces such as electric and magnetic fields (cf. \cite%
{castellanos1998electrohydrodynamics}), vibrational forces of the container
(cf. \cite{Abramson1963}), etc. In those cases, the pressure $p$ at the
interface is not constant, and this results in a nontrivial right hand side
of equation (\ref{le2}), and hence in a nonhomogeneous version of (\ref{le5}%
):
\begin{equation}
\phi _{tt}+\frac{\partial \phi }{\partial n}=h(t,x).  \label{le6}
\end{equation}

\subsection{Conformal transplants}

As is explained in \cite{asmar2002applied,fontelos2020gravity}, let $\Psi$
be a real-valued function written as
\begin{equation*}
\Psi:(x,y)\to\Psi(x,y)=\Psi(z)
\end{equation*}
be defined in a domain $D$. Also, let $\widetilde{\Psi}$ be defined in $%
\widetilde{D}$ as follows: for any $\omega\in\widetilde{D}$ we have
\begin{equation*}
\widetilde{\Psi}(\omega):=\Psi(f^{[-1]}(\omega))=\Psi(x(x^{\prime
},y^{\prime }),y(x^{\prime },y^{\prime })).
\end{equation*}

Then the following relation between the normal derivatives holds:
\begin{equation}  \label{ct8}
\frac{\partial\widetilde{\Psi}}{\partial\widetilde{n}}=\left|\frac{dz}{%
d\omega} \right|\frac{\partial\Psi}{\partial n}=\frac{1}{|f^{\prime }(z)|}%
\frac{ \partial\Psi}{\partial n}.
\end{equation}

\subsection{Reformulation as an integrodifferential PDE}

In the particular case of the half-cylinder geometry, $\Omega$, we use the
conformal map
\begin{equation*}  \label{ie1}
\omega=f(z)=\frac{1}{\frac{1}{2}\left(z+\frac{1}{z}\right)},
\end{equation*}
into the half plane $\Omega^{\prime }:=\{(x^{\prime },y^{\prime })\in\mathbb{%
\ R}^2:y^{\prime }<0\}$.

From (\ref{ct8}), equations (\ref{le1}), (\ref{le4}), and (\ref{le5}), in
variables $\omega=x^{\prime }+iy^{\prime }$ become (see Figure \ref{F1})
\begin{align}
\Delta \widetilde{\phi }& =0,\quad \text{for}\ y^{\prime }<0,  \label{ie5} \\
\widetilde{\phi }_{tt}+|f^{\prime }(x^{\prime })|\frac{\partial \widetilde{%
\phi }}{\partial \widetilde{n}}& =0,\quad \text{at }y^{\prime }=0,\
|x^{\prime }|\leq 1,  \label{ie52} \\
\frac{\partial \widetilde{\phi }}{\partial \widetilde{n}}& =\frac{\widetilde{%
j}(t,x^{\prime })}{|f^{\prime }(x^{\prime })|},\quad \text{at }y^{\prime
}=0,\ |x^{\prime }|>1,  \label{ie53}
\end{align}%
where $j=\widetilde{j}\circ f$. That is, for $|x^{\prime }|>1$, $x^{\prime
}=f(x)$.

By relation (\ref{ct8}) and since $|f^{\prime }(z)|\underset{%
z\to-i}{\longrightarrow}\infty$, we complement system above with the
following boundary conditions at infinity:
\begin{equation*}
\partial_{x^{\prime }}\widetilde{\phi}, \partial_{y^{\prime }}\widetilde{%
\phi }\to 0, \quad \text{as } y^{\prime }\to-\infty \ \text{or } |x^{\prime
}|\to\infty.
\end{equation*}

By taking the Fourier transform in (\ref{ie5}), in the variable $x^{\prime }$; using notation $\widetilde{\Phi }$ and the convention
\begin{equation*}
\Phi {(k)=\frac{1}{\sqrt{2\pi }}\int_{\mathbb{R}}\phi (x^{\prime
})e^{-ikx^{\prime }}dx^{\prime }},
\end{equation*}%
we have
\begin{equation*}
\widetilde{\Phi }_{y^{\prime }y^{\prime }}-k^{2}\widetilde{\Phi }=0,
\end{equation*}%
which implies
\begin{align}
\widetilde{\Phi }(t,k,y^{\prime })& =\widetilde{\Phi }(t,k,0)e^{|k|y^{\prime
}}  \label{ie6} \\
& =\widetilde{\Phi }(t,k,0)\widehat{\frac{-\sqrt{2}y^{\prime }}{\sqrt{\pi }%
(x^{\prime 2}+y^{\prime 2})}}.  \notag
\end{align}

By taking inverse Fourier transform
\begin{equation*}
\widetilde{\phi}(t,x^{\prime },y^{\prime })=-\frac{y^{\prime}}{\pi}
\int_{-\infty}^{+\infty}\frac{\widetilde{\phi}(t,\xi,0)}{(x^{\prime}-%
\xi)^2+y^{\prime 2}}d\xi.
\end{equation*}

Since we want to establish properties on the normal derivative, $\widetilde{%
\phi}_{y^{\prime }}$, taking the $y^{\prime }$ derivative in (\ref{ie6})
instead, and evaluating at $y^{\prime }=0$ we find
\begin{align}
\widetilde{\Phi}_{y^{\prime}}(t,k,0)&=\frac{1}{i}\si(k)(ik)\widetilde{\Phi}
(t,k,0),  \notag \\
&=\widehat{\frac{\sqrt{2}}{\sqrt{\pi}x^{\prime }}}\widetilde{\Phi}
_{x^{\prime }}(t,k,0).  \label{ie7}
\end{align}

Then, taking inverse Fourier transform
\begin{equation}  \label{ie8}
\left.\widetilde{\phi}_{y^{\prime }}\right|_{y^{\prime }=0}=\frac{1}{\pi}
P.V.\int_{-\infty}^{+\infty}\frac{\widetilde{\phi}_{x^{\prime }}(t,\xi,0)}{
x^{\prime}-\xi}d\xi=H\left(\left.\widetilde{\phi}_{x^{\prime
}}\right|_{y^{\prime }=0}\right).
\end{equation}

Since $HH=-I$ (see \cite{hochstadt2011integral}):
\begin{equation}
\left. \widetilde{\phi }_{x^{\prime }}\right\vert _{y^{\prime }=0}=-\frac{1}{%
\pi }P.V.\int_{-\infty }^{\infty }\frac{\widetilde{\phi }_{y^{\prime
}}(t,\xi ,0)}{x^{\prime }-\xi }d\xi =-\frac{1}{\pi }P.V.\int_{-\infty
}^{\infty }\frac{\frac{\partial \widetilde{\phi }}{\partial \widetilde{n}}%
(t,\xi ,0)}{x^{\prime }-\xi }d\xi .  \label{ie9}
\end{equation}

\subsection{The Dirichlet-Neumann operator}

Let us consider in this section the basic situation when $\frac{\partial
\phi }{\partial n}=0$, on the solid walls of the domain $\Omega $. The value
of $\frac{\partial \phi }{\partial n}$ at the fluid interface may be viewed
as the result of an operator (the so-called Dirichlet-Neumann operator)
acting on the function $\phi $ restricted to the interface. We are going to
deduce a few simple consequences obtained from the properties of this
operator. Keeping in mind the equation (\ref{le3}), we can easily deduce the
following mass conservation relation:
\begin{equation}
\frac{d}{dt}\int_{-1}^{1}\zeta dx=\int_{-1}^{1}\frac{\partial \phi }{
\partial n}dx=\int_{\partial \Omega }\frac{\partial \phi }{\partial n}
=\int_{\Omega }\nabla \cdot \left( \nabla \phi \right) =0.  \label{dn1}
\end{equation}%
Notice then that, in the mapped coordinates, by (\ref{ct8})
\begin{equation}
0=\int_{-1}^{1}\frac{\partial \phi }{\partial n}dx=\int_{-1}^{1}\left\vert
f^{\prime }(x^{\prime })\right\vert \frac{\partial \widetilde{\phi }}{
\partial \widetilde{n}}\frac{dx^{\prime }}{\left\vert f^{\prime }(x^{\prime
})\right\vert }=\int_{-1}^{1}\frac{\partial \widetilde{\phi }}{\partial
\widetilde{n}}dx^{\prime },  \label{dn2}
\end{equation}%
and hence, from (\ref{dn1}) and (\ref{ct8})
\begin{equation*}
\frac{d}{dt}\int_{-1}^{1}\widetilde{\zeta }\frac{dx^{\prime }}{\left\vert
f^{\prime }(x^{\prime })\right\vert }=0.
\end{equation*}

Next, since $\frac{\partial \phi }{\partial n}=0$ at the solid boundaries of
$\Omega $, we have
\begin{equation*}
\int_{-1}^{1}\psi \frac{\partial \phi }{\partial n}dx=\int_{\partial \Omega
}\psi \frac{\partial \phi }{\partial n}dx
\end{equation*}%
for any harmonic function $\psi $ also vanishing at the solid boundaries of $%
\Omega $, and therefore by Green's identity we deduce
\begin{equation*}
\int_{\partial \Omega }\psi \frac{\partial \phi }{\partial n}
dx-\int_{\partial \Omega }\phi \frac{\partial \psi }{\partial n}
dx=\int_{\Omega }(\psi \Delta \phi -\phi \Delta \psi)=0,
\end{equation*}%
which implies
\begin{equation*}
\int_{\partial \Omega }\psi \frac{\partial \phi }{\partial n}
dx=\int_{\partial \Omega }\phi \frac{\partial \psi }{\partial n}dx,
\end{equation*}%
or, equivalently,
\begin{equation*}
\int_{-1}^{1}\widetilde{\psi }\frac{\partial \widetilde{\phi }}{\partial
\widetilde{n}}dx^{\prime }=\int_{-1}^{1}\widetilde{\phi }\frac{\partial
\widetilde{\psi }}{\partial \widetilde{n}}dx^{\prime },
\end{equation*}%
showing that the Dirichlet-Neumann operator, mapping $\widetilde{\phi }$ at
the free surface into $\frac{\partial \widetilde{\phi }}{\partial \widetilde{
n}}$, is formally selfadjoint. The details on the
appropriate functional space where the operator is self-adjoint will be
given in section \ref{S4}.

Finally, from (\ref{le2}), (\ref{dn1}), and the mass
conservation property (\ref{mass}):
\begin{equation}  \label{dn4}
\int_{-1}^{1}\phi_{t}dx=-\int_{-1}^{1}\zeta dx=0,
\end{equation}
which implies
\begin{equation*}
\frac{d}{dt}\int_{-1}^{1}\phi(t,x)dx=0.
\end{equation*}
Therefore, by choosing $\phi(0,x)$ such that $\int_{-1}^{1}\phi(0,x)dx=0$
(this can always be achieved by adding a suitable constant to a given $\phi$%
), we have
\begin{equation*}
\int_{-1}^{1}\phi(t,x)dx=0,
\end{equation*}
implying
\begin{equation}  \label{dn5}
\int_{-1}^{1}\widetilde{\phi}(t,x^{\prime })\frac{dx^{\prime }}{\left\vert
f^{\prime }(x^{\prime })\right\vert }=0.
\end{equation}

\section{An associated Cauchy problem}

\label{S3}

As in the last section, let $\frac{\partial \phi }{\partial n}=0$ on the
solid walls of the domain $\Omega $. Then, from the integrodifferential
formulation (\ref{ie9}), we obtain
\begin{equation}
\left. \widetilde{\phi }_{x^{\prime }}\right\vert _{y^{\prime }=0}=-\frac{1}{%
\pi }P.V.\int_{-1}^{1}\frac{\widetilde{\phi }_{y^{\prime }}(t,\xi ,0)}{%
x^{\prime }-\xi }d\xi .  \label{cp1}
\end{equation}%
We observe two properties of $\widetilde{\phi }_{x^{\prime }}$. Firstly, by making use of the identity (see section 4.3 of \cite{tricomi1985integral})
\begin{equation}
\frac{1}{\pi }P.V.\int_{-1}^{1}\frac{1}{\sqrt{1-x^{\prime 2}}}\frac{%
dx^{\prime }}{x^{\prime }-\xi }=0,\ \text{for }\xi \in
\left(-1,1\right),
\end{equation}
we prove%
\begin{align}
\int_{-1}^{1}\frac{\widetilde{\phi }_{x^{\prime }}(x^{\prime })}{\sqrt{%
1-x^{\prime 2}}}dx^{\prime }& =-\frac{1}{\pi }\int_{-1}^{1}\frac{1}{\sqrt{%
1-x^{\prime 2}}}\left( P.V.\int_{-1}^{1}\frac{\widetilde{\phi }_{y^{\prime
}}(\xi )}{x^{\prime }-\xi }d\xi \right) dx^{\prime }  \label{n1} \\
& =-\frac{1}{\pi }\int_{-1}^{1}\widetilde{\phi }_{y^{\prime }}(\xi )\left(
P.V.\int_{-1}^{1}\frac{1}{\sqrt{1-x^{\prime 2}}(x^{\prime }-\xi )}dx^{\prime
}\right) d\xi  \notag \\
& =0.  \notag
\end{align}

Secondly, by the mass conservation $\int_{-1}^{1}%
\widetilde{\phi }_{y^{\prime }}(\xi )d\xi =0 $ (see (\ref{dn2})), and using (see (\ref{a-3}), for $r=1$)
\begin{equation*}
\frac{1}{\pi }P.V.\int_{-1}^{1}\frac{x^{\prime }}{\sqrt{1-x^{\prime 2}}}%
\frac{dx^{\prime }}{x^{\prime }-\xi }=1,\ \text{for }\xi \in
(-1,1),
\end{equation*}%
we also have
\begin{align}
\int_{-1}^{1}\frac{x^{\prime }\widetilde{\phi }_{x^{\prime }}(x^{\prime })}{%
\sqrt{1-x^{\prime 2}}}dx^{\prime }& =-\frac{1}{\pi }\int_{-1}^{1}\frac{%
x^{\prime }}{\sqrt{1-x^{\prime 2}}}\left( P.V.\int_{-1}^{1}\frac{\widetilde{%
\phi }_{y^{\prime }}(\xi )}{x^{\prime }-\xi }d\xi \right) dx^{\prime }
\label{n2} \\
& =-\frac{1}{\pi }\int_{-1}^{1}\widetilde{\phi }_{y^{\prime }}(\xi )\left(
P.V.\int_{-1}^{1}\frac{x^{\prime }}{\sqrt{1-x^{\prime 2}}(x^{\prime }-\xi )}%
dx^{\prime }\right) d\xi  \notag \\
& =-\int_{-1}^{1}\widetilde{\phi }_{y^{\prime }}(\xi )d\xi =0.  \notag
\end{align}

By using the inverse transform instead, from (\ref{cp1}) we have, for $%
|x^{\prime }|<1$,
\begin{equation}
\left. \widetilde{\phi }_{y^{\prime }}\right\vert _{y^{\prime }=0}=\frac{1}{%
\sqrt{1-x^{\prime 2}}}\frac{1}{\pi }P.V.\int_{-1}^{1}\sqrt{1-\xi ^{2}}\frac{%
\widetilde{\phi }_{x^{\prime }}(t,\xi ,0)}{x^{\prime }-\xi }d\xi +\frac{C}{%
\sqrt{1-x^{\prime 2}}},  \label{cp2}
\end{equation}%
where $C$ is an arbitrary constant (see \cite%
{hochstadt2011integral}, Chapter 5.2, Example 13 on the airfoil equation). In order for (\ref{cp2}) to be the general solution to the integral equation
(\ref{cp1}) it suffices to require $\widetilde{\phi }_{x^{\prime
}}(t,x^{\prime },0)$ to satisfy (see {\cite{hochstadt2011integral}})
\begin{equation}
\|\widetilde{\phi }_{x^{\prime }}(t,x^{\prime },0)\|
_{L_{(1-x^{\prime 2})^{1/2}}^{2}}^{2}\equiv \int_{-1}^{1}\sqrt{1-x^{\prime 2}%
}|\widetilde{\phi }_{x^{\prime }}(t,x^{\prime },0)|^{2}dx^{\prime }<\infty \text{.}  \label{norm1}
\end{equation}
Finally, due to mass conservation and using (\ref{n2}), the constant must be
chosen to be zero:
\begin{align*}
0& =\int_{-1}^{1}\widetilde{\phi }_{y^{\prime }}(x^{\prime })dx^{\prime } \\
& =\int_{-1}^{1}\sqrt{1-\xi ^{2}}\widetilde{\phi }_{x^{\prime }}(\xi )\left(
\frac{1}{\pi }P.V.\int_{-1}^{1}\frac{1}{\sqrt{1-x^{\prime 2}}}\frac{1}{%
x^{\prime }-\xi }dx^{\prime }\right) d\xi +\pi C \\
& =\pi C.
\end{align*}%
If, instead of (\ref{norm1}), one assumes the stronger condition
\begin{equation}
\|\widetilde{\phi }_{x^{\prime }}(t,x^{\prime },0)\|
_{L_{(1-x^{\prime 2})^{-1/2}}^{2}}^{2}\equiv \int_{-1}^{1}\frac{|\widetilde{\phi }_{x^{\prime }}(t,x^{\prime },0)|^{2}}{\sqrt{%
1-x^{\prime 2}}}dx^{\prime }<\infty \text{,}  \label{norm2}
\end{equation}%
together with (\ref{n1}), then (cf. {\cite{hochstadt2011integral}, Chapter
5.2, Example 13)}
\begin{equation*}
\left. \widetilde{\phi }_{y^{\prime }}\right\vert _{y^{\prime }=0}=\sqrt{%
1-x^{\prime 2}}\frac{1}{\pi }P.V.\int_{-1}^{1}\frac{1}{\sqrt{1-\xi ^{2}}}%
\frac{\widetilde{\phi }_{x^{\prime }}(t,\xi ,0)}{x^{\prime }-\xi }d\xi .
\end{equation*}

As a consequence of the computations above, it is possible to
state the following.
\begin{lemma}
Given $\widetilde{\phi }_{x^{\prime }}\in L_{\sqrt{1-x^{\prime 2}}}^{2}(-1,1)$, let $\widetilde{\phi }_{y^{\prime }}$ be satisfying (\ref{cp1}) and the mass
conservation condition $\int_{-1}^{1}\widetilde{\phi }_{y^{\prime }}(\xi
)d\xi =0$. Then, the following relations hold,
\begin{align}
\left. \widetilde{\phi }_{y^{\prime }}\right\vert _{y^{\prime }=0}& =\frac{1%
}{\sqrt{1-x^{\prime 2}}}\frac{1}{\pi }P.V.\int_{-1}^{1}\sqrt{1-\xi ^{2}}%
\frac{\widetilde{\phi }_{x^{\prime }}(t,\xi ,0)}{x^{\prime }-\xi }d\xi
\label{ex2} \\
& =\partial _{x^{\prime }}\left( \sqrt{1-x^{\prime 2}}\frac{1}{\pi }%
P.V.\int_{-1}^{1}\frac{1}{\sqrt{1-\xi ^{2}}}\frac{\widetilde{\phi }(t,\xi ,0)%
}{x^{\prime }-\xi }d\xi \right) .  \label{ex3}
\end{align}%
Moreover, if $\widetilde{\phi }_{x^{\prime }}\in L_{(1-x^{\prime
2})^{-1/2}}^{2}(-1,1)\subset L_{\sqrt{1-x^{\prime 2}}}^{2}(-1,1)$ then $%
\widetilde{\phi }_{y^{\prime }}$ may also be given by the equivalent expression
\begin{equation}
\left. \widetilde{\phi }_{y^{\prime }}\right\vert _{y^{\prime }=0}=\sqrt{%
1-x^{\prime 2}}\frac{1}{\pi }P.V.\int_{-1}^{1}\frac{1}{\sqrt{1-\xi ^{2}}}%
\frac{\widetilde{\phi }_{x^{\prime }}(t,\xi ,0)}{x^{\prime }-\xi }d\xi.
\label{ex1}
\end{equation}
\end{lemma}

\begin{proof}
We present a proof based on expansions in terms of Tchebyshev polynomials in
the Appendix \ref{app1}.
\end{proof}

Expression (\ref{ex1}) was used in our previous article \cite%
{fontelos2020gravity}, while expression (\ref{ex3}) will be more convenient
in the present work. We present next an alternative deduction of (\ref{ex3}%
). We introduce the function $\chi$ defined by
\begin{equation*}
\widetilde{\phi }=\chi _{x^{\prime }}
\end{equation*}%
so that
\begin{equation*}
\chi =\int_{-\infty }^{x^{\prime }}\widetilde{\phi }(s,y^{\prime })ds,
\end{equation*}%
(note that integrability at $-\infty $ is guarantied by the decay of $\widetilde{\phi}$ from (\ref{ie9}) after integration by parts and the mass conservation (\ref{dn2})) and
\begin{equation*}
\Delta \chi =0.
\end{equation*}

We have then
\begin{equation*}
\chi _{y^{\prime }}=\frac{\partial }{\partial y^{\prime }}\int_{-\infty
}^{x^{\prime }}\widetilde{\phi }dx^{\prime }=\int_{-\infty }^{x^{\prime
}}(\chi _{y^{\prime }})_{x^{\prime }}dx^{\prime }=
\begin{cases}
0, & x^{\prime }\leq -1, \\
\chi _{y^{\prime }}, & -1<x^{\prime }\leq 1, \\
0, & x^{\prime }>1.%
\end{cases}%
\end{equation*}%
Since $\Delta \chi =0$, by (\ref{cp1}) we have
\begin{equation}
\chi _{x^{\prime }}(x^{\prime })=-\frac{1}{\pi }P.V.\int_{-1}^{1}\frac{%
\chi_{y^{\prime }}(\xi )}{x^{\prime }-\xi }d\xi .  \label{f1}
\end{equation}

Therefore, inverting Hilbert transform
\begin{equation}
\chi _{y^{\prime }}(x^{\prime })=\sqrt{1-x^{\prime 2}}\frac{1}{\pi }
P.V.\int_{-1}^{1}\frac{1}{\sqrt{1-\xi ^{2}}}\frac{\chi _{x^{\prime }}(\xi )}{
x^{\prime }-\xi }d\xi ,  \label{f2}
\end{equation}%
and taking $x^{\prime }$ derivative
\begin{equation}
\widetilde{\phi }_{y^{\prime }}(x^{\prime })=\partial _{x^{\prime }}\left(
\sqrt{1-x^{\prime 2}}\frac{1}{\pi }P.V.\int_{-1}^{1}\frac{1}{\sqrt{1-\xi
^{2} }}\frac{\widetilde{\phi }(\xi )}{x^{\prime }-\xi }d\xi \right)
\label{f3}
\end{equation}%
which, by (\ref{ie52}), yields the following evolution problem
\begin{equation}
\widetilde{\phi }_{tt}=-\left\vert f^{\prime }(x^{\prime })\right\vert
\partial _{x^{\prime }}\left( \sqrt{1-x^{\prime 2}}\frac{1}{\pi }
P.V.\int_{-1}^{1}\frac{\widetilde{\phi }(t,\xi )}{\sqrt{1-\xi ^{2}}
(x^{\prime }-\xi )}d\xi \right) .  \label{integro}
\end{equation}

The equation (\ref{integro}) may be rewritten as an identical equation for $%
\widetilde{\zeta }$ after taking an additional time derivative and writing $%
\widetilde{\phi }_{t}=-\widetilde{\zeta }$. We remark again that the
formulation (\ref{integro}) is equivalent to the formulation given in \cite%
{fontelos2020gravity} (see Appendix \ref{app1}).

We need to complement (\ref{integro}) with suitable initial and boundary
conditions, namely
\begin{align*}
\widetilde{\phi}(0,x^{\prime })&=\widetilde{\phi}_{0}(x^{\prime }), \\
\widetilde{\phi}_{t}(0,x^{\prime })&=\widetilde{\phi}_{1}(x^{\prime }),
\end{align*}
where, by equation (\ref{dn4})
\begin{equation}  \label{condi1}
\int_{-1}^{1}\frac{\widetilde{\phi}_{1}(x^{\prime })}{\left\vert
f^{\prime}(x^{\prime })\right\vert}dx^{\prime }={\color{blue}-}\int_{-1}^{1}%
\frac{ \widetilde{\zeta}(x^{\prime })}{\left\vert f^{\prime }(x^{\prime
})\right\vert }dx^{\prime }=0.
\end{equation}

Moreover, from (\ref{n1}), (\ref{n2}), we impose on the initial data:
\begin{equation}
\int_{-1}^{1}\frac{\widetilde{\phi }_{0,x^{\prime }}(x^{\prime })}{\sqrt{
1-x^{\prime 2}}}dx^{\prime }=0,  \label{condi2}
\end{equation}
\begin{equation}
\int_{-1}^{1}\frac{x^{\prime }\widetilde{\phi }_{0,x^{\prime }}(x^{\prime })
}{\sqrt{1-x^{\prime 2}}}dx^{\prime }=0,  \label{condi2'}
\end{equation}%
which are automatically fulfilled by defining an initial normal derivative $%
\widetilde{\phi }_{0,y^{\prime }}(x^{\prime })$ with vanishing mean value in
$\left[-1,1\right] $ and the corresponding tangential derivative $\widetilde{%
\phi }_{0,x^{\prime }}(x^{\prime })$ defined by (\ref{cp1}).

In the case when $\left\vert f^{\prime}(x^{\prime })\right\vert$ is a
symmetric function (corresponding to a symmetric domain $\Omega$) it is
useful to think of $\widetilde{\phi}(t,x^{\prime })$ as decomposed into
symmetric and antisymmetric part; that is
\begin{equation*}
\widetilde{\phi}(t,x^{\prime })=S(t,x^{\prime })+N(t,x^{\prime })
\end{equation*}
with
\begin{equation*}
S(t,x)=\frac{\widetilde{\phi}(t,x)+\widetilde{\phi}(t,-x)}{2}, \quad N(t,x)=
\frac{\widetilde{\phi}(t,x)-\widetilde{\phi}(t,-x)}{2},
\end{equation*}
and the initial data decomposed accordingly
\begin{align*}
S_{0}(x)&=\frac{\widetilde{\phi}_{0}(x)+\widetilde{\phi}_{0}(-x)}{2},\quad
N_{0}(x)=\frac{\widetilde{\phi}_{0}(x)-\widetilde{\phi}_{0}(-x)}{2}, \\
S_{1}(x)&=\frac{\widetilde{\phi}_{1}(x)+\widetilde{\phi}_{1}(-x)}{2}, \quad
N_{1}(x)=\frac{\widetilde{\phi}_{1}(x)-\widetilde{\phi}_{1}(-x)}{2}.
\end{align*}

Then, one can consider the evolution problem for symmetric and antisymmetric
functions separately. For symmetric functions, the initial data need to
satisfy (\ref{condi1}), (\ref{condi2'}) while for antisymmetric functions
the conditions are (\ref{condi1}), (\ref{condi2}). Finally, we remark that (%
\ref{condi1}), (\ref{condi2}), (\ref{condi2'}) do not only hold initially,
but for any time by replacing $(\widetilde{\phi}_{0}(x^{\prime }),\widetilde{%
\phi}_{1}(x^{\prime }))$ by $(\widetilde{\phi}(t,x^{\prime }),\widetilde{\phi%
}_{t}(t,x^{\prime }))$.

We discuss now on boundary conditions for (\ref{integro}). Two kind of
conditions are customary: pinned end and free end boundary conditions. In
pinned end conditions one imposes $\widetilde{\zeta}(t,\pm 1)=0$, and since $%
\widetilde{\phi}_{t}=-\widetilde{\zeta}$, this implies
\begin{equation}  \label{pinned}
\widetilde{\phi}(t,\pm 1)=0\ \text{(pinned-end boundary condition)}.
\end{equation}
In the case of free-end boundary conditions one imposes $\widetilde{\zeta}%
_{x^{\prime }}(t,\pm 1)=0$ and since $\widetilde{\phi}_{x^{\prime }t}=%
\widetilde{\zeta}_{x^{\prime }}$ this translates into the condition
\begin{equation}  \label{free}
\widetilde{\phi}_{x^{\prime }}(t,\pm 1)=0\ \ \text{(free-end boundary
condition)}.
\end{equation}

We discuss now the case with fluid injection, i.e. with the condition (\ref%
{ie53}) where the flux $j$ is such that$\ \int \frac{\widetilde{j}%
(t,x^{\prime })}{|f^{\prime }(x^{\prime })|}dx^{\prime }=0$ (in order to
preserve the total fluid mass). Then, by (\ref{ie9}), equation (\ref{cp1})
needs to be replaced by%
\begin{equation*}
\left. \widetilde{\phi }_{x^{\prime }}\right\vert _{y^{\prime }=0}=-\frac{1}{%
\pi }P.V.\int_{-1}^{1}\frac{\widetilde{\phi }_{y^{\prime }}(t,\xi ,0)}{%
x^{\prime }-\xi }d\xi +\frac{1}{\pi }\int_{\mathbb{R}\setminus \lbrack -1,1]}%
\frac{\widetilde{j}(t,z)}{|f^{\prime }(z)|(\xi -z)}dz.
\end{equation*}%
Inverting as in (\ref{ex2}), (\ref{ex3}) we get the following formula for
the normal derivative:
\begin{align}
\left. \widetilde{\phi }_{y^{\prime }}\right\vert _{y^{\prime }=0}&
=\partial _{x^{\prime }}\left( \sqrt{1-x^{\prime 2}}\frac{1}{\pi }%
P.V.\int_{-1}^{1}\frac{1}{\sqrt{1-\xi ^{2}}}\frac{\widetilde{\phi }(t,\xi ,0)%
}{x^{\prime }-\xi }d\xi \right)  \notag \\
& \phantom{= }-\frac{1}{\sqrt{1-x^{\prime 2}}}\frac{1}{\pi ^{2}}%
P.V.\int_{-1}^{1}\frac{\sqrt{1-\xi ^{2}}}{(x^{\prime }-\xi )}\int_{\mathbb{R}%
\setminus \lbrack -1,1]}\frac{\widetilde{j}(t,z)}{|f^{\prime }(z)|(\xi -z)}%
dzd\xi ,
\end{align}%
which leads to the evolution equation
\begin{equation}
\widetilde{\phi }_{tt}=-\left\vert f^{\prime }(x^{\prime })\right\vert
\partial _{x^{\prime }}\left( \sqrt{1-x^{\prime 2}}\frac{1}{\pi }%
P.V.\int_{-1}^{1}\frac{\widetilde{\phi }(t,\xi )}{\sqrt{1-\xi ^{2}}%
(x^{\prime }-\xi )}d\xi \right) +\widetilde{h}(t,x^{\prime }),
\label{integro2}
\end{equation}%
where
\begin{equation}
\widetilde{h}(t,x^{\prime })=\frac{\left\vert f^{\prime }(x^{\prime
})\right\vert }{\sqrt{1-x^{\prime 2}}}\frac{1}{\pi ^{2}}P.V.\int_{-1}^{1}%
\frac{\sqrt{1-\xi ^{2}}}{(x^{\prime }-\xi )}\int_{\mathbb{R}\setminus
\lbrack -1,1]}\frac{\widetilde{j}(t,z)}{|f^{\prime }(z)|(\xi -z)}dzd\xi .
\label{hxt}
\end{equation}

Notice that $\chi _{y^{\prime }}=\int_{-\infty }^{x^{\prime }}\widetilde{%
\phi }_{y^{\prime }}dx^{\prime }$ implies expression (\ref{hxt}) can be
obtained by replacing (\ref{f1}) with
\begin{equation*}
\chi _{x^{\prime }}(x^{\prime })=-\frac{1}{\pi }P.V.\int_{-1}^{1}\frac{\chi
_{y^{\prime }}(\xi )}{x^{\prime }-\xi }d\xi -\frac{1}{\pi } P.V.\int_{%
\mathbb{R}\setminus \lbrack -1,1]}\frac{J(\xi )}{x^{\prime }-\xi }d\xi
\end{equation*}%
where $J$ is the primitive of $\widetilde{j}/|f^{\prime }|$, inverting the
Hilbert transform (restricted to $\left[ -1,1\right] $) in the first term at
the right hand side as in (\ref{f2}), (\ref{f3}) and using finally (\ref{ex2})--(\ref{ex3}).

As mentioned above, equation (\ref{integro2}), for given $\widetilde{h}%
(t,x^{\prime })$, is also valid in situations where the interface is
actuated by means of external forces such as electric and magnetic ones,
external container vibration, etc. Hence, we will present a general
discussion on controllability for general $\widetilde{h}(t,x^{\prime })$ and
will only specify for boundary injection in the final section.

\section{Spectrum of the sloshing problem}

\label{S4}

From now on, from equation (\ref{integro2}) and to summarize
the computations from the preceding sections, we are concerned with the
following Initial Value Problem,
\begin{equation}
\begin{cases}
\frac{\widetilde{\phi }_{tt}}{|f^{\prime }|}+\mathcal{A}\widetilde{\phi }=
\widetilde{h}(t,x^{\prime }),\quad & (t,x^{\prime })\in (0,\infty )\times
(-1,1), \\
\widetilde{\phi }(0,x^{\prime })=\widetilde{\phi }_{0}(x^{\prime }), &
x^{\prime }\in (-1,1), \\
\widetilde{\phi }_{t}(0,x^{\prime })=\widetilde{\phi }_{1}(x^{\prime }), &
x^{\prime }\in (-1,1),%
\end{cases}
\label{cp}
\end{equation}
where $\mathcal{A}$ is the integral, non-local operator
\begin{equation*}
\mathcal{A}\widetilde{\phi}\equiv \partial_{x^{\prime }}\left(\sqrt{
1-x^{\prime 2}}\frac{1}{\pi}P.V.\int_{-1}^{1}\frac{\widetilde{\phi}(\xi)}{
\sqrt{1-\xi^{2}}(x^{\prime }-\xi)}d\xi\right).
\end{equation*}

Following \cite{micu2004introduction}, to study the interior controllability problem (\ref{cp}), we need to consider the homogeneous (backward in time)
adjoint version in $(0,T)$ as that given in (\ref{integro}). For that purpose, let us
consider first the eigenvalue problem
\begin{equation}  \label{eigen}
\lambda\frac{\widetilde{\phi}}{\left\vert f^{\prime}(x^{\prime })\right\vert}
=\partial_{x^{\prime }}\left(\sqrt{1-x^{\prime 2}}P.V.\frac{1}{\pi}
\int_{-1}^{1}\frac{\widetilde{\phi}(\xi)}{\sqrt{1-\xi^{2}}(x^{\prime }-\xi)}
d\xi\right).
\end{equation}

Let $T_n(x)$, $U_n(x)$, with $n\in\mathbb{N}\cup\{0\}$, be the Tchebyshev polynomials of the first and second kind respectively. They are
defined as the polynomial solutions of the equations (see
\cite{abramowitz1972} for details)
\begin{align*}
T_n(\cos\theta)&=\cos(n\theta), \\
U_n(\cos\theta)&=\frac{\sin((n+1)\theta)}{\sin\theta}.
\end{align*}
They satisfy the following orthogonality relations in $L^2(-1,1)$ with the
corresponding inner products $\langle f,g\rangle=\int_{-1}^1\frac{fg}{\sqrt{%
1-x^2}}dx$, $\langle f,g\rangle=\int_{-1}^1fg\sqrt{1-x^2}dx$:
\begin{equation}\label{oT}
\int_{-1}^1T_n(x)T_m(x)\frac{dx}{\sqrt{1-x^2}}=
\begin{cases}
0, \ \text{if }n\ne m, \\
\pi \ \text{if }n=m=0, \\
\frac{\pi}{2} \ \text{if }n=m\ne 0,%
\end{cases}%
\end{equation}
\begin{equation}\label{oU}
\int_{-1}^1U_n(x)U_m(x)\sqrt{1-x^2}dx=
\begin{cases}
0, \ \text{if }n\ne m, \\
\frac{\pi}{2} \ \text{if }n=m.%
\end{cases}%
\end{equation}
Moreover, for $n\ge1$, we have relations
\begin{align}
\frac{1}{\pi}P.V.\int_{-1}^{1}\frac{T_{n}(\xi)}{\sqrt{1-\xi^{2}}(x^{\prime
}-\xi)}d\xi&=-U_{n-1}(x^{\prime }),  \label{a-3} \\
\frac{d}{dx^{\prime }}\left(\sqrt{1-x^{\prime 2}}U_{n-1}(x^{\prime
})\right)&=-r\frac{T_{n}(x^{\prime })}{\sqrt{1-x^{\prime 2}}}.  \notag
\end{align}

By writing
\begin{equation}  \label{a-1}
\widetilde{\phi}(\xi)=\sum_{n=0}^{\infty}a_{n}T_{n}(\xi),
\end{equation}
and using the identities above, we find
\begin{equation}  \label{aphi}
\partial_{x^{\prime }}\left(\sqrt{ 1-x^{\prime 2}}\frac{1}{\pi}%
P.V.\int_{-1}^{1}\frac{\widetilde{\phi}(\xi)}{ \sqrt{1-\xi^{2}}(x^{\prime
}-\xi)}d\xi\right) =\sum_{n=1}^{\infty}na_{n} \frac{T_{n}(x^{\prime })}{%
\sqrt{1-x^{\prime 2}}}.
\end{equation}
We use now the fact that $\int_{-1}^{1}\widetilde{\phi}(x^{\prime })\frac{
dx^{\prime }}{\left\vert f^{\prime}(x^{\prime })\right\vert}=0$ (see (\ref%
{dn5})), to deduce
\begin{equation*}
a_{0}\int_{-1}^{1}T_{0}(x^{\prime })\frac{dx^{\prime }}{\left\vert
f^{\prime}(x^{\prime })\right\vert}+\sum_{n=1}^{\infty}a_{n}
\int_{-1}^{1}T_{n}(x^{\prime })\frac{dx^{\prime }}{\left\vert
f^{\prime}(x^{\prime })\right\vert}=0.
\end{equation*}
So that
\begin{equation*}
a_{0}=-\frac{1}{\int_{-1}^{1}T_{0}(x^{\prime })\frac{dx^{\prime }}{
\left\vert f^{\prime}(x^{\prime })\right\vert}}\sum_{n=1}^{\infty}a_{n}
\int_{-1}^{1}T_{n}(x^{\prime })\frac{dx^{\prime }}{\left\vert
f^{\prime}(x^{\prime })\right\vert},
\end{equation*}
and conclude the estimate%
\begin{equation*}
a_{0}^{2}\leq\frac{\sum_{n=1}^{\infty}\left(\int_{-1}^{1}T_{n}(x^{\prime })
\frac{dx^{\prime }}{\left\vert f^{\prime}(x^{\prime })\right\vert }
\right)^{2}}{\left(\int_{-1}^{1}T_{0}(x^{\prime })\frac{dx^{\prime }}{
\left\vert f^{\prime}(x^{\prime })\right\vert}\right)^{2}}
\sum_{n=1}^{\infty}a_{n}^{2}.
\end{equation*}

Now, let $c_n$ be such that
\begin{equation*}
\frac{\sqrt{1-x^{\prime 2}}}{|f^{\prime }(x^{\prime })|}=\sum_{n=0}^{
\infty}c_nT_n(x^{\prime }).
\end{equation*}
Then, by the orthogonality of the Tchebyshev polynomials,
\begin{equation*}
c_0=\frac{1}{\pi}\int_{-1}^{1}\frac{T_0(x^{\prime })}{|f^{\prime }(x^{\prime
})|}dx^{\prime }, \quad c_n=\frac{2}{\pi}\int_{-1}^{1}\frac{T_n(x^{\prime })
}{|f^{\prime }(x^{\prime })|}dx^{\prime }.
\end{equation*}

Therefore
\begin{align*}
\int_{-1}^1\frac{1-x^{\prime 2}}{|f^{\prime}(x^{\prime})|^2}\frac{1}{\sqrt{
1-x^{\prime 2}}}dx^{\prime }=\sum_{n,m=0}^{\infty}c_nc_m\int_{-1}^{1}\frac{
T_nT_m}{\sqrt{1-x^{\prime 2}}}dx^{\prime }=\pi c_0^2+\frac{\pi}{2}
\sum_{n=1}^{\infty}c_n^2,
\end{align*}
namely
\begin{equation*}
\frac{1}{\pi}\left(\int_{-1}^{1}T_{0}(x^{\prime })\frac{dx^{\prime }}{
\left\vert f^{\prime}(x^{\prime })\right\vert}\right)^{2}+\frac{2}{\pi}
\sum_{n=1}^{\infty}\left(\int_{-1}^{1}T_{n}(x^{\prime })\frac{dx^{\prime }}{
\left\vert f^{\prime}(x^{\prime })\right\vert }\right)^{2} =\int_{-1}^{1}
\sqrt{1-x^{\prime 2}}\frac{dx^{\prime }}{\left\vert f^{\prime}(x^{\prime
})\right\vert^{2}}<\infty,
\end{equation*}
which implies
\begin{equation}  \label{a0}
a_{0}^{2}\leq C\sum_{n=1}^{\infty}a_{n}^{2}.
\end{equation}

Let
\begin{equation}  \label{psi3}
\widetilde{\psi}=\sum_{n=0}^\infty b_nT_n(x^{\prime }).
\end{equation}
We consider now the scalar product
\begin{equation*}
(\mathcal{A}\widetilde{\phi},\widetilde{\psi})_{L^{2}}=\int_{-1}^{1}
\widetilde{\psi}\frac{\partial\widetilde{\phi} }{\partial \widetilde{n}}
dx^{\prime },
\end{equation*}
and find
\begin{equation}  \label{a1}
(\mathcal{A}\widetilde{\phi},\widetilde{\psi})_{L^{2}}=\int_{-1}^{1}\left(
\sum_{m=0}^{\infty}b_{m}T_{m}(x^{\prime })\right)
\left(\sum_{n=1}^{\infty}na_{n}\frac{T_{n}(x^{\prime })}{\sqrt{1-x^{\prime 2}%
}}\right)dx^{\prime }=\frac{\pi}{2}\sum_{n=1}^{\infty }na_{n}b_{n},
\end{equation}
so that, $(\widetilde{\psi},\mathcal{A}\widetilde{\phi}) _{L^{2}}=(\mathcal{A%
}\widetilde{\psi},\widetilde{\phi}) _{L^{2}}$, implying the selfadjoint
character of the operator $\mathcal{A}$, a fact already shown in a previous
section.

We study the problem
\begin{equation}  \label{phiu}
\mathcal{A}\widetilde{\phi}=u,
\end{equation}
where $u\in L_{w}^{2}$ with $w=\sqrt{1-x^{\prime 2}}$ and
\begin{equation*}
L_{w}^{2}\equiv \left\{ u:\int_{-1}^{1}\sqrt{1-x^{\prime 2}}\left\vert
u\right\vert ^{2}dx^{\prime }<\infty\right\}.
\end{equation*}
Since $\left\{T_{n}(x^{\prime }){/\sqrt{1-x^{\prime 2}}}\right\}
_{n=0}^{\infty}$ form an orthogonal base for $L_{w}^{2}$, given $u\in L_w^2$
we can expand
\begin{equation*}
u=\sum_{n=0}^{\infty}u_{n}\frac{T_{n}(x^{\prime })}{\sqrt{1-x^{\prime 2}}},
\end{equation*}
and hence, assuming
\begin{equation}
u_{0}=\int_{-1}^{1}u(x^{\prime })dx^{\prime }=0,  \label{u0null}
\end{equation}
write the following weak version of (\ref{phiu})
\begin{equation}  \label{weak1}
(\mathcal{A}\widetilde{\phi},\widetilde{\psi})_{L^{2}}=(u,\widetilde{\psi}
)_{L^{2}},
\end{equation}
to be satisfied for any $\widetilde{\psi}$ in
\begin{equation}  \label{psi1}
H_{w^{-1}}^{\frac{1}{2}}\equiv \left\{\widetilde{\psi} :\left\Vert\widetilde{%
\psi}\right\Vert_{H_{w^{-1}}^{\frac{1}{2}}}^2=\left\Vert \widetilde{\psi}%
\right\Vert_{L_{w^{-1}}^2}^2+\sum_{n=1}^{\infty}n\left( \int_{-1}^{1}%
\frac{T_{n}(x^{\prime })}{\sqrt{1-x^{\prime 2}}}\widetilde{\psi} (x^{\prime
})dx^{\prime }\right)^{2}<\infty \right\}
\end{equation}
such that, in addition,
\begin{equation}  \label{psi2}
\int_{-1}^{1}\widetilde{\psi }(x^{\prime })\frac{dx^{\prime }}{\left\vert
f^{\prime }(x^{\prime })\right\vert }=0.
\end{equation}

By using the preceding functional framework, we can state a result on
existence of weak solutions for problem (\ref{phiu}).

\begin{lemma}
Let $u\in L^2_{w}$ satisfying (\ref{u0null}). Then, there exists a unique
weak solution $\widetilde{\phi}\in H_{w^{-1}}^{1/2}$ for problem (\ref{phiu}%
).
\end{lemma}

\begin{proof}
We can write the equation (\ref{weak1}) in the form
\begin{equation}  \label{weak2}
\sum_{n=1}^{\infty }na_{n}b_{n}=\sum_{n=1}^{\infty }u_{n}b_{n},
\end{equation}
making it clear that $(\mathcal{A}\widetilde{\phi},\widetilde{\psi})_{L^{2}}$
defines a continuous bilinear form on $H_{w^{-1}}^{\frac{1}{2}}$. Indeed,

\begin{remark}
By using the orthogonality property (\ref{oT}) on (\ref{psi3}) we have, for $n\ge1$,
\begin{equation*}
\sqrt{n}\int_{-1}^{1}\frac{T_{n}(x^{\prime })}{\sqrt{1-x^{\prime 2}}}
\widetilde{\psi}(x^{\prime })dx^{\prime }=\frac{\pi}{2}\sqrt{n}b_n,
\end{equation*}
or, equivalently,
\begin{equation*}
\sum_{n=1}^{\infty}n\left(\int_{-1}^{1}\frac{T_{n}(x^{\prime })}{\sqrt{
1-x^{\prime 2}}}\widetilde{\psi}(x^{\prime })dx^{\prime
}\right)^{2}=\frac{\pi^2}{4}\sum_{n=1}^{\infty}n|b_n|^2.
\end{equation*}

Notice, from this last equality and (\ref{a1}), that $%
H_{w^{-1}}^{1/2}$ can be written as
\begin{equation*}
H_{w^{-1}}^{\frac{1}{2}}\equiv \left\{\widetilde{\psi} :\left\Vert\widetilde{
\psi}\right\Vert_{H_{w^{-1}}^{\frac{1}{2}}}^2=\left\Vert \widetilde{\psi}
\right\Vert_{L_{w^{-1}}^2}^2+\frac{\pi}{2}(\mathcal{A}\widetilde{\psi},\widetilde{\psi})_{L^2}<\infty \right\}.
\end{equation*}
\end{remark}

Thus from (\ref{a1}) and the last remark, we have
\begin{align*}
|(\mathcal{A}\widetilde{\phi},\widetilde{\psi})_{L^2}|&=\frac{\pi}{2}
\left|\sum_{n\ge1}na_nb_n\right| \le\frac{\pi}{2}\left(\sum_{n\ge1}n|a_n|^2
\right)^{1/2}\left(\sum_{n\ge1}n|b_n|^2\right)^{1/2} \\
&\le \frac{2}{\pi}\|\widetilde{\phi}\|_{H_{w^{-1}}^{1/2}}\|\widetilde{\psi}
\|_{H_{w^{-1}}^{1/2}}.
\end{align*}

$(\mathcal{A}\widetilde{\phi},\widetilde{\psi})_{L^{2}}$ is also coercive
(the $L_{w^{-1}}^{2}$ part of the norm of $\widetilde{ \phi}$ is trivially
bounded by $\sum_{n=1}^{\infty }na_{n}^{2}$, except for $a_{0}^{2}$ that is
also bounded by (\ref{a0})):

\begin{remark}
By using the orthogonality property (\ref{oT}) on (\ref{a-1}),
\begin{equation*}
\|\widetilde{\phi}\|_{L_{w^{-1}}^2}^2=\int_{-1}^{1}\frac{|\widetilde{\phi}%
|^2 }{\sqrt{1-x^{\prime 2}}}dx^{\prime }=\pi a_0^2+\frac{\pi}{2}
\sum_{n=1}^{\infty}a_n^2.
\end{equation*}
Moreover,
\begin{equation*}
a_n=
\begin{cases}
\frac{1}{\pi}\int_{-1}^{1}\frac{T_n(x^{\prime })}{\sqrt{1-x^{\prime 2}}}
\widetilde{\phi}(x^{\prime })dx^{\prime }, \quad & n=0, \\
\frac{2}{\pi}\int_{-1}^{1}\frac{T_n(x^{\prime })}{\sqrt{1-x^{\prime 2}}}
\widetilde{\phi}(x^{\prime })dx^{\prime }, & n\ge1.%
\end{cases}%
\end{equation*}
\end{remark}

Then, from the last remark and (\ref{a0}), coerciveness holds as follows:
\begin{align}  \label{a8}
(\mathcal{A}\widetilde{\phi},\widetilde{\phi})_{L^2}&=\frac{\pi}{2}
\sum_{n\ge1}n|a_n|^2 \ge C\left[\pi\sum_{n=1}^{\infty}a_n^2+\frac{\pi}{2}
\sum_{n=1}^{\infty}a_n^2+\frac{\pi^2}{4}\sum_{n=1}^{\infty}na_n^2\right] \\
&\ge C\left[\pi a_0^2+\frac{\pi}{2}\sum_{n=1}^{\infty}a_n^2+\frac{\pi^2}{4}
\sum_{n=1}^{\infty}na_n^2\right]  \notag \\
&=C\left[\|\widetilde{\phi}\|_{L_{w^{-1}}^2}^2+\frac{\pi}{2}(\mathcal{A}\widetilde{\phi},\widetilde{\phi})\right]  \notag \\
&=C\|\widetilde{\phi}\|_{H_{w^{-1}}^{1/2}}^2,  \notag
\end{align}
for a constant $C>0$.

Finally, $u\in L_{w}^{2}$ defines a linear continuous functional on $%
H_{w^{-1}}^{1/2}$:
\begin{align*}
\left|\int_{-1}^{1}u\widetilde{\phi}dx^{\prime }\right|&\le \int_{-1}^{1}
\sqrt[4]{1-x^{\prime 2}}|u|\frac{\widetilde{\phi}}{\sqrt[4]{1-x^{\prime 2}}}
dx^{\prime } \\
&\le\|u\|_{L_w^2}\|\widetilde{\phi}\|_{L_{w^{-1}}^2} \\
&\le\|u\|_{L_w^2}\|\widetilde{\phi}\|_{H_{w^{-1}}^{1/2}}.
\end{align*}

Hence, by Lax-Milgram's theorem, there exists a unique weak solution $%
\widetilde{\phi}$ to (\ref{phiu}) that belongs to $H_{w^{-1}}^{\frac{1}{2}}$.
\end{proof}

If we consider now $u=\frac{v}{\left\vert f^{\prime}(x^{\prime })\right\vert
}$ then the problem
\begin{equation*}
\mathcal{A}\widetilde{\phi}=\frac{v}{\left\vert f^{\prime }(x^{\prime
})\right\vert }
\end{equation*}
where $v$ is such that
\begin{equation}
\int_{-1}^{1}\sqrt{1-x^{\prime 2}}\frac{v^{2}(x^{\prime })}{\left\vert
f^{\prime}(x^{\prime })\right\vert^{2}}dx^{\prime }<\infty  \label{cond1}
\end{equation}
has also a unique solution $\widetilde{\phi}\in H_{w^{-1}}^{\frac{1}{2}}$.
Note that (\ref{cond1}) is satisfied if $v\in L_{\left\vert
f^{\prime}(x^{\prime })\right\vert ^{-1}}^{2}$ provided there exists a
constant $C$ such that
\begin{equation}  \label{cond2}
\left\vert f^{\prime}(x^{\prime })\right\vert \geq C\sqrt{1-x^{\prime 2}},
\quad x^{\prime }\in[-1,1],
\end{equation}
because
\begin{equation*}
\int_{-1}^{1}\sqrt{1-x^{\prime 2}}\frac{v^{2}(x^{\prime })}{\left\vert
f^{\prime}(x^{\prime })\right\vert^{2}}dx^{\prime }\leq\frac{1}{C^{2}}\int_{-1}^{1}\frac{
v^{2}(x^{\prime })}{\sqrt{1-x^{\prime 2}}}dx^{\prime },
\end{equation*}
or
\begin{equation*}
\left\Vert u\right\Vert_{L_{|f^{\prime}|^{-1}}^{2}}^{2}\leq \frac{1}{C}
\left\Vert u\right\Vert_{L_{w^{-1}}^{2}}^{2}.
\end{equation*}

We define now the operator $\mathcal{T}$ such that $\widetilde{\phi}=%
\mathcal{T}v$. If we view $\mathcal{T}$ as an operator from $L_{\left\vert
f^{\prime }\right\vert ^{-1}}^{2}$ to $L_{\left\vert f^{\prime}\right\vert
^{-1}}^{2}$, since $H_{w^{-1}}^{\frac{1}{2}}\subset L_{w^{-1}}^{2}\subset
L_{\left\vert f^{\prime}\right\vert^{-1}}^{2}$, we can see that $\mathcal{T}$
is a selfadjoint operator (due to the selfadjoint character of $\mathcal{A}$%
) and also a compact operator (due to the compact embedding $H_{w^{-1}}^{%
\frac{1}{2}}\subset L_{w^{-1}}^{2}$ (see Appendix \ref{app2}) and the
continuous inclusion $L_{w^{-1}}^{2}\subset L_{\left\vert
f^{\prime}\right\vert ^{-1}}^{2}$) . Moreover, $\ker (\mathcal{T})=\left\{
0\right\}$ and
\begin{equation*}
(\mathcal{T}v,v)_{L_{|f^{\prime}|^{-1}}^2}=(\widetilde{\phi},|f^{\prime }|
\mathcal{A}\widetilde{\phi})_{L_{|f^{\prime}|^{-1}}^2}=(\widetilde{\phi},
\mathcal{A}\widetilde{\phi})_{L^2}\ge0,
\end{equation*}
for any $v\in L_{\left\vert f^{\prime}\right\vert ^{-1}}^{2}$. Therefore, by
the spectral decomposition theorem, $L_{\left\vert f^{\prime }\right\vert
^{-1}}^{2}$ admits a Hilbert basis $\left\{ e_{n}\right\} $ formed by
eigenvectors of $\mathcal{T}$, with eigenvalues $\mu _{n}$ such that $\mu
_{n}>0$ and $\mu _{n}\rightarrow 0$ as $n\rightarrow\infty$. We have then {$%
\mathcal{T}e_n=\mu_ne_n$}, $e_{n}\in H_{w^{-1}}^{\frac{1}{2}}$, and thus
\begin{equation}  \label{a10}
{(\mathcal{A}e_{n},\widetilde{\psi})_{L^{2}}=\frac{1}{\mu_n}(\mathcal{A}(
\mathcal{T}e_{n}),\widetilde{\psi})_{L^{2}}=\frac{1}{\mu_{n}}\left(\frac{
e_{n}}{|f^{\prime }|},\widetilde{\psi}\right)_{L^{2}},} \quad \text{for any }
\widetilde{\psi}\in H_{w^{-1}}^{\frac{1}{2}}.
\end{equation}

This implies $e_{n}$ is a weak solution to the eigenvalue problem (\ref%
{eigen}) with $\lambda =\lambda _{n}=\mu _{n}^{-1}$. The base can be made
orthonormal, i.e.
\begin{equation*}
\int_{-1}^{1}\frac{e_{i}e_{j}}{\left\vert f^{\prime}(x^{\prime })\right\vert}
dx^{\prime }=\delta _{ij}.
\end{equation*}

Let us remark that condition (\ref{cond2}) is satisfied provided the angles $%
\alpha _{1,2}$ between the fluid interface and the solid container satisfy
\begin{equation*}
\alpha _{1,2}\geq \frac{\pi }{2},
\end{equation*}
since $|f^{\prime }(x^{\prime })|$ has a singularity weaker (i.e. with
larger exponent) than $|x\pm1|^{1/2}$ otherwise \cite{asmar2002applied}.

We summarize the results above in the following theorem. {\color{blue} }

\begin{theorem}
Let $\Omega $ be a domain such that the interior angles between the free
liquid interface and the solid are greater or equal than $\pi /2$, so that (\ref{cond2}) is fulfilled. There exist a Hilbert basis
$\{e_{n}\}_{n\geq 1}$ of $L_{\left\vert f^{\prime}\right\vert ^{-1}}^{2}$
and a sequence $\{\lambda _{n}\}_{n\geq 1}$ of real numbers with $\lambda
_{n}>0 $ $\forall n$ and $\lambda _{n}\rightarrow +\infty $ such that
\begin{equation*}
e_{n}\in H_{w^{-1}}^{\frac{1}{2}},
\end{equation*}
\begin{equation*}
\mathcal{A}e_{n}=\lambda _{n}\frac{e_{n}}{\left\vert f^{\prime }(x^{\prime
})\right\vert }.
\end{equation*}
\end{theorem}

We will show next that the eigenfunctions $e_{n}$ possess
further regularity. From (\ref{a-1})-(\ref{aphi}), if $e_{n}(x^{\prime
})=\sum_{m=0}^{\infty }a_{m}T_{m}(x^{\prime })$, by using (\ref{oT}) we have
\begin{equation*}
\int_{-1}^{1}\sqrt{1-x^{\prime 2}}\left\vert \mathcal{A}e_{n}\right\vert
^{2}dx^{\prime }=\frac{\pi }{2}\sum_{m=1}^{\infty }m^{2}a_{m}^{2}.
\end{equation*}

Therefore, by using condition (\ref{cond2}),
\begin{align*}
\int_{-1}^{1}\sqrt{1-x^{\prime 2}}\left\vert \mathcal{A}e_{n}\right\vert
^{2}dx^{\prime } &=\lambda _{n}^{2}\int_{-1}^{1}\sqrt{1-x^{\prime 2}}\frac{
\left\vert e_{n}\right\vert ^{2}}{\left\vert f^{\prime }(x^{\prime
})\right\vert ^{2}}dx^{\prime } \\
&\leq C_{n}\int_{-1}^{1}\frac{\left\vert e_{n}\right\vert ^{2}}{\sqrt{
1-x^{\prime 2}}}dx^{\prime }=C_{n}\Vert e_{n}\Vert _{L_{w^{-1}}^{2}}^{2},
\end{align*}%
we have that $\sum_{m=1}^{\infty }m^{2}a_{m}^{2}$ is bounded by $%
\|e_n\|_{L^2_{w^{-1}}}^2$. On the other hand, since $\frac{dT_{m}(x^{\prime
})}{dx^{\prime }}=mU_{m-1}(x^{\prime })$, from the orthogonality property (\ref{oU}) we deduce%
\begin{equation*}
\int_{-1}^{1}\sqrt{1-x^{\prime 2}}\left\vert \frac{de_{n}}{dx'}\right\vert
^{2}dx^{\prime }=\frac{\pi }{2}\sum_{m=1}^{\infty }m^{2}a_{m}^{2},
\end{equation*}%
which is bounded by $\Vert e_{n}\Vert_{L_{w^{-1}}^{2}}^{2}$. Hence, $e_{n}\in H_{w^{-1}}^{1}$ implying, by Sobolev embeddings, that $e_{n}(x^{\prime
}) $ is a continuous function in $(-1,1)$. Moreover, it is bounded at $x=\pm1$ since:
\begin{equation*}
e_{n}(x^{\prime })-e_{n}(-1)=\int_{-1}^{x^{\prime
}}e_{n}^{\prime }(x)dx\leq \left( \int_{-1}^{x^{\prime }}\frac{1}{\sqrt{
1-x^{2}}}dx\right) ^{\frac{1}{2}}\left( \int_{-1}^{x^{\prime }}\sqrt{1-x^{2}}
\left\vert e_{n}^{\prime }(x)\right\vert ^{2}dx\right) ^{\frac{1}{2}}.
\end{equation*}

So that
\begin{equation*}
\sup_{-1\leq x^{\prime }\leq -1+\delta }\frac{\left\vert
e_{n}(x')-e_{n}(-1)\right\vert ^{2}}{\left\vert 1+x'\right\vert ^{\frac{1}{2}}}
\leq C\int_{-1}^{1}\sqrt{1-x^2}\left\vert e_{n}^{\prime
}\right\vert ^{2}dx,
\end{equation*}%
and hence%
\begin{equation*}
\left\vert e_{n}(x^{\prime })-e_{n}(-1)\right\vert \leq C\left\vert
1+x'\right\vert ^{\frac{1}{4}}.
\end{equation*}

Identical estimate may be obtained in the neighborhood of $x'=1$.

Note that no condition has been imposed on the eigenfunctions at $x'=\pm 1$.
The eigenfunctions correspond to free-end boundary conditions. We can also
consider the pinned-end boundary condition by writing the same weak
formulation (\ref{weak1}) but assuming that $\widetilde{\psi }$ belongs to
the closure of $C_{c}^{\infty }(-1,1)$ in the topology defined by the $%
H_{w^{-1}}^{\frac{1}{2}}$ norm and satisfying (\ref{psi2}). We denote such
space as $H_{w^{-1},0}^{\frac{1}{2}}$. The same arguments as for the
free-end case lead to the existence of a complete set of eigenfunctions as
in the Theorem above.

As a final remark, notice that for bounded symmetric $\left\vert
f^{\prime}(x^{\prime })\right\vert$, one has $L_{\left\vert f^{\prime
}(x^{\prime })\right\vert^{-1}}^{2}\subset L^{2}$ and hence vectors in $%
L_{\left\vert f^{\prime}\right\vert ^{-1}}^{2}$ may be expanded in
trigonometric basis of $L^{2}$. Using the base $\left\{ \sin(n\pi x^{\prime
})\right\} _{n=1}^{\infty }$ one can approximate antisymmetric
eigenfunctions (solutions to (\ref{eigen})) with pinned contact lines while
the set $\left\{\cos\left(\left( n+\frac{1}{2}\right)\pi x^{\prime }\right)
\right\}_{n=0}^{\infty}$ (with the extra mass conservation condition also
imposed) allows to approximate symmetric eigenfunctions. Likewise, the sets $%
\left\{\sin\left(\left(n+\frac{1}{2}\right)\pi x^{\prime }\right)
\right\}_{n=0}^{\infty}$ and $\left\{\cos(n\pi x^{\prime
})\right\}_{n=1}^{\infty}$ allow to find antisymmetric and symmetric
eigenfunctions with free-end condition. This approach was followed in \cite%
{fontelos2020gravity}, in order to compute eigenvalues and eigenfunctions
numerically for various domains.

\section{Observability and energy estimates}

\label{S5}

We will consider next the control problem consisting of finding $\widetilde{h}(t,x^{\prime})$ in equation (\ref{integro2}) such that for given initial
data the solution vanishes at some time $T$. As it is customary in
controllability theory, we need to consider first the homogeneous adjoint
problem and obtain observability estimates.

The analysis presented above allows us to consider the solution in $(0,T)$
of the following homogeneous adjoint problem (see (\ref{integro}))
\begin{equation}
\frac{\widetilde{\phi }_{tt}}{\left\vert f^{\prime }\right\vert }+\mathcal{A}%
\widetilde{\phi }=0,  \label{o1}
\end{equation}%
with suitable initial conditions $\widetilde{\phi }(0)=\widetilde{\phi }_{0}$%
, $\widetilde{\phi }^{\prime }(0)=\widetilde{\phi }_{1}$, whose solution we
can write as
\begin{equation}
\widetilde{\phi }=\sum_{n=1}^{\infty }\left( A_{n}\cos (\theta
_{n}t)+B_{n}\sin (\theta _{n}t)\right) e_{n},  \label{o2}
\end{equation}%
where $\theta _{n}=\sqrt{\lambda _{n}}$ and
\begin{equation*}
\widetilde{\phi }_{0}=\sum_{n=1}^{\infty }A_{n}e_{n},\quad \widetilde{\phi }%
_{1}=\sum_{n=1}^{\infty }\theta _{n}B_{n}e_{n},
\end{equation*}%
with
\begin{equation*}
A_{n}=\int_{-1}^{1}\frac{\widetilde{\phi }_{0}e_{n}}{\left\vert f^{\prime
}(x^{\prime })\right\vert }dx^{\prime }=(\widetilde{\phi }%
_{0},e_{n})_{|f^{\prime}|^{-1}},\quad B_{n}=\frac{1}{\theta _{n}}%
\int_{-1}^{1}\frac{\widetilde{\phi }_{1}e_{n}}{\left\vert f^{\prime
}(x^{\prime })\right\vert }dx^{\prime }=\frac{1}{\theta _{n}}(\widetilde{%
\phi }_{1},e_{n})_{|f^{\prime}|^{-1}}.
\end{equation*}

Hence $\Vert \widetilde{\phi }\Vert _{L^{2}(L_{|f^{\prime
}|^{-1}}^{2})}^{2}=\int_{0}^{T}\int_{-1}^{1}\frac{\widetilde{\phi }^{2}}{%
\left\vert f^{\prime }(x^{\prime })\right\vert }dx^{\prime }dt$ can be
bounded from below by
\begin{align}
& \int_{0}^{T}\sum_{n}\left( A_{n}^{2}\cos ^{2}(\theta _{n}t)+B_{n}^{2}\sin
^{2}(\theta _{n}t)+2A_{n}B_{n}\sin (\theta _{n}t)\cos (\theta _{n}t)\right)
dt  \label{o6} \\
& =\sum_{n}\left[ \frac{A_{n}^{2}}{2}\left( T+\frac{\sin (2\theta _{n}T)}{%
2\theta _{n}}\right) +\frac{B_{n}^{2}}{2}\left( T-\frac{\sin (2\theta _{n}T)%
}{2\theta _{n}}\right) +A_{n}B_{n}\frac{1-\cos (2\theta _{n}T)}{2\theta _{n}}%
\right]  \notag \\
& \geq \sum_{n}\left[ \frac{A_{n}^{2}}{2}\left( T+\frac{\sin (2\theta _{n}T)%
}{2\theta _{n}}\right) +\frac{B_{n}^{2}}{2}\left( T-\frac{\sin (2\theta
_{n}T)}{2\theta _{n}}\right) -\frac{A_{n}^{2}+B_{n}^{2}}{2}\frac{1-\cos
(2\theta _{n}T)}{2\theta _{n}}\right]  \notag \\
& =\frac{T}{2}\sum_{n}\left[ A_{n}^{2}\left( 1-\frac{1}{2\theta _{n}T}+\frac{%
\sin (2\theta _{n}T)}{2\theta _{n}T}+\frac{\cos (2\theta _{n}T)}{2\theta
_{n}T}\right) \right.  \notag \\
& \phantom{\ =}\left. +B_{n}^{2}\left( 1-\frac{1}{2\theta _{n}T}-\frac{\sin
(2\theta _{n}T)}{2\theta _{n}T}+\frac{\cos (2\theta _{n}T)}{2\theta _{n}T}%
\right) \right]  \notag \\
& \geq C\frac{T}{2}\sum_{n}\left[ A_{n}^{2}+B_{n}^{2}\right] ,  \notag
\end{align}%
for some positive $C>0$ and $T>2.42/2\min (\theta _{n})=2.42/2\theta _{1}$.

From the orthogonality of the Hilbert basis $\{e_n\}$,
\begin{equation}  \label{o7}
\|\widetilde{\phi}_0\|_{L_{|f^{\prime}|^{-1}}^2}^2=\int_{-1}^{1}\frac{
\widetilde{\phi}_0^2}{|f^{\prime }(x^{\prime })|}dx^{\prime
}=\sum_{n=1}^{\infty}A_n^2.
\end{equation}
Moreover, from (\ref{a10}), the decomposition given in (\ref{o2}), and (\ref%
{o1}) we conclude $\lambda_{n}=\theta_n^2$ and thus
\begin{equation}  \label{o8}
\sum_{n=1}^{\infty}B_{n}^{2}=\sum_{n}\left(\frac{1}{\theta_{n}}\int_{-1}^{1}
\frac{\widetilde{\phi}_{1}e_{n}}{\left\vert f^{\prime}(x^{\prime
})\right\vert }dx^{\prime }\right)^{2} =\sum_{n}\frac{1}{\lambda _{n}}
\left(\int_{-1}^{1}\frac{\widetilde{\phi}_{1}e_{n}}{\left\vert
f^{\prime}(x^{\prime })\right\vert}dx^{\prime }\right)^{2}.
\end{equation}

Let us define, accordingly, the space
\begin{equation*}
\widetilde{H}^{-\frac{1}{2}}=\left\{\widetilde{\phi}:\sum_{n}\frac{1}{
\lambda_{n}}\left(\int_{-1}^{1}\frac{\widetilde{\phi}e_{n}}{\left\vert
f^{\prime }(x^{\prime })\right\vert}dx^{\prime }\right)^{2}<\infty\right\},
\end{equation*}
and its dual space
\begin{equation*}
\widetilde{H}^{\frac{1}{2}}=\left\{\widetilde{\phi}:\sum_{n}\lambda_{n}
\left(\int_{-1}^{1}\frac{\widetilde{\phi}e_{n}}{\left\vert f^{\prime
}(x^{\prime })\right\vert}dx^{\prime }\right)^{2}<\infty\right\}.
\end{equation*}
Note that, by (\ref{a1}) and definition of $\lambda _{n}$ and $e_{n}$,
\begin{equation*}
{\ \left\Vert\widetilde{\phi}\right\Vert _{\widetilde{H}^{\frac{1}{2}}}^{2}
=\sum \left( \int \frac{\widetilde{\phi} \lambda _{n}e_{n}}{\left\vert
f^{\prime}(x)\right\vert }\right) \left( \int \frac{\widetilde{\phi} e_{n}}{
\left\vert f^{\prime }(x)\right\vert }\right) =\int \widetilde{\phi}\mathcal{%
\ A}\widetilde{\phi} =\frac{\pi}{2}\sum_{n\ge1} na_{n}^{2} \sim\left\Vert
\widetilde{\phi} \right\Vert _{H_{w^{-1}}^{\frac{1}{2}}}^{2}. }
\end{equation*}
From (\ref{o6})--(\ref{o8}), we conclude then for some $C>0$ and $T$
sufficiently large,
\begin{equation*}
C^{-1}T\left( \left\Vert\widetilde{\phi}_{1}\right\Vert_{\widetilde{H}^{-%
\frac{1}{2}}}^{2}+\left\Vert\widetilde{\phi}_{0}\right\Vert_{L_{|f^{%
\prime}|^{-1}}^2}^{2}\right) \ge \int_{0}^{T}\int_{-1}^{1}\frac{\widetilde{%
\phi}^{2}}{\left\vert f^{\prime }(x^{\prime })\right\vert }dx^{\prime }dt
\ge CT\left(\left\Vert\widetilde{\phi}_{1}\right\Vert_{\widetilde{H}^{-\frac{%
1}{2}}}^{2}+\left\Vert\widetilde{\phi}_{0}\right\Vert_{L_{|f^{%
\prime}|^{-1}}^2}^{2}\right).
\end{equation*}

This proves the following theorem:

\begin{theorem}
The equation (\ref{o1}) is observable in time $T>T_0$ for some $T_0$
sufficiently large. That is, there exists a positive constant $C>0$ such
that for $T$ sufficiently large
\begin{equation*}
\int_{0}^{T}\int_{-1}^{1}\frac{\widetilde{\phi}^{2}}{\left\vert f^{\prime
}(x^{\prime })\right\vert }dx^{\prime }dt\ge CT\left( \left\Vert\widetilde{
\phi}_{1}\right\Vert_{\widetilde{H}^{-\frac{1}{2}}}^{2}+\left\Vert\widetilde{
\phi}_{0}\right\Vert_{L_{|f^{\prime}|^{-1}}^2}^{2}\right).
\end{equation*}
\end{theorem}

We consider next the nonhomogeneous problem (as in (\ref{integro2}))
\begin{equation}  \label{o12}
\frac{\widetilde{\psi}_{tt}}{\left\vert f^{\prime}\right\vert}+\mathcal{A}
\widetilde{\psi}=\frac{\widetilde{h}}{|f^{\prime }|}.
\end{equation}
It will be convenient the following energy estimate:
\begin{equation}  \label{o13}
\frac{1}{2}\frac{d}{dt}\int_{-1}^{1}\frac{\widetilde{\psi}_{t}^{2}}{
\left\vert f^{\prime }(x^{\prime })\right\vert}dx^{\prime }+\int_{-1}^{1}
\widetilde{\psi}_{t}\mathcal{A}\widetilde{\psi}dx^{\prime }=\int_{-1}^{1}
\frac{\widetilde{h}\widetilde{\psi}_{t}}{|f^{\prime }(x^{\prime })|}
dx^{\prime }.
\end{equation}

If we let
\begin{equation*}
\widetilde{\psi}=\sum_{n=1}^{\infty}c_ne_n,
\end{equation*}
then
\begin{align*}
c_{n}=\int_{-1}^{1}&\frac{\widetilde{\psi}e_{n}}{\left\vert f^{\prime
}(x^{\prime })\right\vert}dx^{\prime }, \\
\widetilde{\psi}_t=\sum_{n=1}^{\infty}c_{n,t}e_n \quad\implies\quad
&\int_{-1}^{1}\frac{\psi_{t}^{2}}{\left\vert f^{\prime }(x^{\prime
})\right\vert}dx^{\prime }=\sum_{n=1}^{\infty}c_{n,t}^{2}.
\end{align*}
Noticing that
\begin{align*}
\int_{-1}^{1}\widetilde{\psi}_{t}\mathcal{A}\widetilde{\psi}%
dx^{\prime}&=\int_{-1}^{1}\left(\sum_{n}c_{n,t}e_n\right)\left(\sum_{m}%
\lambda_mc_m\frac{e_m}{|f^{\prime }|}\right)dx^{\prime} \\
&=\sum_{n}\lambda_{n}c_{n,t}c_{n} \\
&=\frac{1}{2}\frac{d}{dt}\sum_{n}\lambda_{n}c_{n}^{2},
\end{align*}
(\ref{o13}) is equivalent to the following inequality
\begin{equation*}
\frac{dE}{dt}\le\|\widetilde{h}\|_{L_{|f^{\prime}|^{-1}}^2}E^{\frac{1}{2}},
\end{equation*}
for the energy defined as
\begin{equation*}
E=\frac{1}{2}\sum_{n}c_{n,t}^{2}+\frac{1}{2}\sum_{n}\lambda _{n}c_{n}^{2}.
\end{equation*}

Hence, the natural initial data for which the problem is well-posed is $%
\left(\widetilde{\psi}_{0},\widetilde{\psi}_{1}\right) \in \widetilde{H}^{
\frac{1}{2}}\times L_{|f^{\prime}|^{-1}}^2$ and, the energy is bounded
provided $\int_{0}^{T}\|\widetilde{h}\|_{L_{|f^{\prime}|^{-1}}^2}^{2}dt<%
\infty $, i.e. $\widetilde{h}\in L^2(0,T;L_{|f^{\prime}|^{-1}}^2)$. Namely,
we have the following: {\color{blue} }

\begin{theorem}
Given $f$ a conformal mapping that satisfies condition (\ref%
{cond2}), for any $\widetilde{h}\in L^2(0,T;L_{|f^{\prime}|^{-1}}^2)$ and $%
\left( \widetilde{\psi}_{0},\widetilde{\psi}_{1}\right) \in \widetilde{H}^{%
\frac{1}{ 2}}\times L_{|f^{\prime}|^{-1}}^2$ equation (\ref{o12}) has a
unique weak solution
\begin{equation*}
(\widetilde{\psi},\widetilde{\psi}^{\prime })\in C([0,T];H_{w^{-1}}^{\frac{ 1%
}{2}}\times L_{|f^{\prime}|^{-1}}^2).
\end{equation*}
\end{theorem}

\section{Controllability}

\label{S6}

Once we have studied in the previous sections the forward evolution equation
and the backward homogeneous system, given by (\ref{o12}) and (\ref{o1})
respectively, we are in position of comparing them to get the
controllability condition on the system.

That is, if we assume the control drives the initial data of system (\ref%
{o12}) to zero, by multiplying the source term in (\ref{o12}) by the
solution to the adjoint problem (\ref{o1}), we obtain
\begin{align*}
\int_{0}^{T}\int_{-1}^{1}\frac{\widetilde{h}\widetilde{\phi}}{|f^{\prime}|}%
dx^{\prime}dt&=\int_{0}^{T}\int_{-1}^{1}\left[\frac{\widetilde{\psi}_{tt}}{%
|f^{\prime}|}+\mathcal{A}\widetilde{\psi}\right]\widetilde{\phi}dx^{\prime
}dt \\
&=\int_{-1}^{1}\frac{\widetilde{\psi}_{t}\widetilde{\phi}|_{0}^{T}}{%
|f^{\prime}|}dx^{\prime}-\int_{0}^{T}\int_{-1}^{1}\frac{\widetilde{\psi} _{t}%
\widetilde{\phi}_{t}}{|f^{\prime }|}dx^{\prime }dt+\int_{0}^{T}\int_{-1}^{1}%
\widetilde{\psi}\mathcal{A}\widetilde{\phi}dx^{\prime}dt \\
&=\int_{-1}^{1}\frac{\widetilde{\psi}_{t}\widetilde{\phi}|_{0}^{T}}{%
|f^{\prime}|}dx^{\prime}-\int_{-1}^{1}\frac{\widetilde{\psi}\widetilde{\phi}%
_{t}|_{0}^{T}}{|f^{\prime }|}dx^{\prime }+\int_{0}^{T}\int_{-1}^{1}
\widetilde{\psi}\left[\frac{\widetilde{\phi}_{tt}}{|f^{\prime }|}+\mathcal{A}
\widetilde{\phi}\right]dx^{\prime }dt \\
&=-\int_{-1}^{1}\frac{\widetilde{\psi}_{1}\widetilde{\phi}_{0}}{|f^{\prime
}| }dx^{\prime }+\int_{-1}^{1}\frac{\widetilde{\psi}_{0}\widetilde{\phi}_{1}%
}{ |f^{\prime }|}dx^{\prime },
\end{align*}
where $(\widetilde{\psi}_{0},\widetilde{\psi}_{1})=(\widetilde{\psi}(0),%
\widetilde{\psi}_{t}(0))\in\widetilde{H}^{\frac{1}{2}}\times
L_{|f^{\prime}|^{-1}}^2$. So that $\widetilde{h}$ can be chosen as the
minimizer of the functional
\begin{equation*}
J[\widetilde{\phi}_{0},\widetilde{\phi}_{1}] =\frac{1}{2}\int_{0}^{T}
\int_{-1}^{1}\frac{\widetilde{\phi}^{2}}{|f^{\prime }|}dx^{\prime
}dt+\int_{-1}^{1}\frac{\widetilde{\psi}_{1}\widetilde{\phi}_{0}}{|f^{\prime
}|}dx^{\prime }-\int_{-1}^{1}\frac{\widetilde{\psi}_{0}\widetilde{\phi}_{1}}{
|f^{\prime }|}dx^{\prime }.
\end{equation*}

Notice that
\begin{equation*}
\left|\int_{-1}^{1}\frac{\widetilde{\psi}_{1}\widetilde{\phi}_{0}}{
|f^{\prime }|}dx^{\prime }-\int_{-1}^{1}\frac{\widetilde{\psi}_{0}\widetilde{
\phi}_{1}}{|f^{\prime }|}dx^{\prime }\right| \le \|\widetilde{\phi}
_{0}\|_{L_{|f^{\prime}|^{-1}}^2}\|\widetilde{\psi}_{1}\|_{L_{|f^{
\prime}|^{-1}}^2}+\|\widetilde{\phi}_{1}\|_{\widetilde{H}^{-\frac{1}{2}}}\|
\widetilde{\psi}_{0}\|_{\widetilde{H}^{\frac{1}{2}}}.
\end{equation*}

As it is customary in controllability theory, coerciveness of the functional $J $ is guaranteed if the adjoint problem is observable in time $T$, that is:
\begin{equation}  \label{obs}
\int_{0}^{T}dt\int_{-1}^{1}\frac{\widetilde{\phi}^{2}}{|f^{\prime }|}
dx^{\prime }\ge C\left(\|\widetilde{\phi}(0)\|_{L_{|f^{\prime}|^{-1}}^2}^2+%
\| \widetilde{\phi}_{t}(0)\|_{\widetilde{H}^{-\frac{1}{2}}}^{2}\right),
\end{equation}
a fact that was proved in the previous section. Hence, we have proved the
following controllability theorem: {\color{blue} }

\begin{theorem}
\label{t4} The system (\ref{o12}) is exactly controllable in time $T$. That
is, for any initial data $\left(\widetilde{\psi}_{0},\widetilde{\psi}%
_{1}\right) \in \widetilde{H}^{\frac{1}{2}}\times L_{|f^{\prime}|^{-1}}^2$,
there exist $\widetilde{h}\in L^2(0,T;L^2_{|f^{\prime}|^{-1}})$ and $T>0$
such that
\begin{equation*}
\|\widetilde{\psi}\|_{H_{w^{-1}}^{\frac{1}{2}}}=0, \quad\text{for }t>T.
\end{equation*}
\end{theorem}

The fact that the control $\widetilde{h}\in
L^{2}((0,T);L_{|f^{\prime}|^{-1}}^2) $ implies that one can write
\begin{equation*}
\widetilde{h}(t,x^{\prime })=\sum_{n}\widetilde{h}_{n}(t)e_{n},
\end{equation*}
with
\begin{equation*}
\int_{0}^{T}\sum_{n}|\widetilde{h}_{n}(t)|^2dt<\infty.
\end{equation*}

We are going to discuss next on how to approach the control function $%
\widetilde{h}$, defined at the fluid interface by means of a function $j$
defined at the solid boundaries. This function $j$ will represent a fluid
injection at certain points $\left\{ x_{1}^{\prime },x_{2}^{\prime
},...,x_{N}^{\prime }\right\} $ of the solid boundary with a flow rate $%
\left\{ j_{1}(t),j_{2}(t),...,j_{N}(t)\right\} $ and under the mass
conservation condition
\begin{equation*}
\sum_{i=1}^{N}j_{i}(t)=0.
\end{equation*}

More precisely,
\begin{equation}
\frac{\partial \widetilde{\phi }}{\partial \widetilde{n}}=\sum_{i}j_{i}(t)
\delta (x^{\prime }-x_{i}^{\prime }).  \label{flux}
\end{equation}
We replace the expression for $\widetilde{j}(t,z,0)$ in (\ref{hxt}) by the
right hand side of (\ref{flux}) at the solid side walls and obtain
\begin{align}
\frac{|f^{\prime }(x^{\prime })|}{\sqrt{1-x^{\prime 2}}}\frac{1}{\pi^{2}}%
&P.V.\int_{-1}^{1}\frac{\sqrt{1-\xi ^{2}}}{(x^{\prime }-\xi )}\int_{\mathbb{%
\ R}\setminus \lbrack -1,1]}\frac{\widetilde{j}(t,z,0)}{|f^{\prime }(z)|(\xi
-z)}dzd\xi  \notag  \label{control} \\
& =\sum_{i}\frac{j_{i}(t)}{|f^{\prime }(x_{i}^{\prime })|}\frac{|f^{\prime
}(x^{\prime })|}{\sqrt{1-x^{\prime 2}}}\frac{1}{\pi ^{2}}P.V.\int_{-1}^{1}
\frac{\sqrt{1-\xi ^{2}}}{(x^{\prime }-\xi )(\xi -x_{i}^{\prime })}d\xi
\notag \\
& =\sum_{i}\frac{j_{i}(t)}{|f^{\prime }(x_{i}^{\prime })|}\frac{|f^{\prime
}(x^{\prime })|}{\sqrt{1-x^{\prime 2}}}\frac{1}{x^{\prime }-x_{i}^{\prime }}
\frac{1}{\pi ^{2}}P.V.\int_{-1}^{1}\sqrt{1-\xi ^{2}}\left[ \frac{1}{
x^{\prime }-\xi }+\frac{1}{\xi -x_{i}^{\prime }}\right] d\xi  \notag \\
& =\frac{|f^{\prime }(x^{\prime })|}{\sqrt{1-x^{\prime 2}}}\sum_{i}j_{i}(t)%
\frac{G(x_{i}^{\prime })-G(x)}{\pi ^{2}|f^{\prime }(x_{i}^{\prime })|}\frac{1%
}{x^{\prime }-x_{i}^{\prime }},
\end{align}%
with
\begin{equation*}
G(x_{i}^{\prime })=P.V.\int_{-1}^{1}\frac{\sqrt{1-\xi ^{2}}}{(\xi
-x_{i}^{\prime })}d\xi .
\end{equation*}

Hence, the challenge is to approximate the control function $\widetilde{h}$
by (\ref{control}). One possibility is to approximate the first $N-1$ modes
so that
\begin{equation*}
\int_{-1}^{1}\frac{\widetilde{h}e_{i}}{\left\vert f^{\prime }(x^{\prime
})\right\vert }dx^{\prime }=h_{i}(t),\quad i=1,...,N-1,
\end{equation*}%
or, denoting $m_{ij}=\frac{1}{\pi ^{2}|f^{\prime }(x_{j}^{\prime })|}%
\int_{-1}^{1}\frac{G(x_{j}^{\prime })-G(x)}{\sqrt{1-x^{\prime 2}}(x^{\prime
}-x_{j}^{\prime })}e_{i}dx^{\prime }$ and including the mass conservation
condition, solving the system%
\begin{align*}
\sum_{j=1}^{N}m_{ij}j_{j}(t)& =h_{i}(t),\quad i=1,...,N-1, \\
\sum_{i=1}^{N}j_{i}(t)& =0.
\end{align*}

\bigskip

\bigskip

The fact that $h_{i}(t)\in L^{2}\left( 0,T\right) $ implies $j_{i}(t)\in
L^{2}\left( 0,T\right) $.

Let us summarize what we have achieved so far in connection
with the injection of fluid problem and the control of splashing appearing
in a cooper converter \cite{brimacombe1985toward,godoy2008modeling}. First,
the main relation between the inner source $\widetilde{h}$, and the
injection of fluids $\widetilde{j}$ at the solid boundary of the half-plane,
is given by (\ref{hxt}). Second, by the computations above, we found the
source $\widetilde{h}$ given in (\ref{flux}), which corresponds to the fluid
injection of jets $\widetilde{j}$ on a finite number of points. Third, once $%
\widetilde{j}$ is known on the half-plane, we can recover $j$, on the
cylindrical container, from (\ref{ct8}). Of course this procedure is
directly connected to the controllability, through Theorem \ref{t4}.

Finally, it is important to highlight that, although we used the case of the
cylinder as a reference, all the computations are valid for any simply
connected domain through the conformal mapping term $|f^{\prime }|$, which
appears as a factor on the main operator for the Cauchy problem.

As far as the three-dimensional case is concerned, the conformal mapping
strategy cannot be replicated directly. While it is true that the integral
equations can be formulated in an equivalent frame in the 3d case (see \cite%
{kim2015capillary} for more details).

\section*{Acknowledgments}

We would like to thank the anonymous referees for their comments and suggestion which helped to significantly improve this work.

\appendix

\section{Equivalent expressions of the operator $\mathcal{A}$}

\label{app1}

We discuss in this appendix the relation between different expressions of
the operator $\mathcal{A}$ based on expansions in Tchebyshev polynomials. We
remind that the Tchebyshev polynomials $\left\{ T_{n}(x^{\prime })\right\}
_{n=0}^{\infty }$ form an orthogonal basis in $L_{(1-x^{\prime
2})^{-1/2}}^{2}(-1,1)$ and the Tchebyshev polynomials $\left\{
U_{n}(x^{\prime })\right\} _{n=0}^{\infty }$ form an orthogonal basis in $%
L_{(1-x^{\prime 2})^{1/2}}^{2}(-1,1)$. \ Remind that tangential and normal
derivatives of $\widetilde{\phi }$ in the interval $(-1,1)$ are related by%
\begin{equation}
\widetilde{\phi }_{x^{\prime }}(x^{\prime })=-\frac{1}{\pi }P.V.\int_{-1}^{1}%
\frac{\widetilde{\phi }_{y^{\prime }}(\xi )}{x^{\prime }-\xi }d\xi .
\label{i0}
\end{equation}

Let
\begin{equation}
\widetilde{\phi }=\sum_{n=0}^{\infty }a_{n}T_{n},  \label{expan}
\end{equation}%
and observe the following well-known identities
\begin{align}
\frac{dT_{n}}{dx^{\prime }}& =nU_{n-1},  \label{i1} \\
-U_{n-1}(x^{\prime })& =\frac{1}{\pi }P.V.\int_{-1}^{1}\frac{T_{n}(\xi )}{%
\sqrt{1-\xi ^{2}}(x^{\prime }-\xi )}d\xi ,  \label{i2} \\
T_{n}(x^{\prime })& =\frac{1}{\pi }P.V.\int_{-1}^{1}\frac{\sqrt{1-\xi ^{2}}%
U_{n-1}(\xi )}{(x^{\prime }-\xi )}d\xi ,  \label{i3} \\
\frac{d}{dx^{\prime }}\left( \sqrt{1-x^{\prime 2}}U_{r-1}(x^{\prime
})\right) & =-r\frac{T_{r}(x^{\prime })}{\sqrt{1-x^{\prime 2}}}.  \label{i4}
\end{align}%
Since, by (\ref{i1}),%
\begin{equation}
\widetilde{\phi }_{x^{\prime }}=\sum_{n=1}^{\infty }na_{n}U_{n-1}(x^{\prime
}),  \label{expphix}
\end{equation}%
we have%
\begin{equation*}
\int_{-1}^{1}\sqrt{1-x^{\prime 2}}|\widetilde{\phi }_{x^{\prime
}}(x^{\prime })|^{2}dx^{\prime }=\frac{\pi }{2}\sum_{n=1}^{\infty}n^{2}a_{n}^{2}.
\end{equation*}

By writing
\begin{equation*}
\widetilde{\phi }_{y^{\prime }}=\sum_{n=0}^{\infty }b_{n}\frac{%
T_{n}(x^{\prime })}{\sqrt{1-x^{\prime 2}}},
\end{equation*}%
and using
\begin{equation*}
\int_{-1}^{1}\widetilde{\phi }_{y^{\prime }}(x^{\prime })dx^{\prime
}=\sum_{n=0}^{\infty }b_{n}\int \frac{T_{0}(x^{\prime })T_{n}(x^{\prime })}{%
\sqrt{1-x^{\prime 2}}}=b_{0}=0,
\end{equation*}%
together with (\ref{i0}), (\ref{i2}) we conclude $b_{n}=na_{n}$. By using (\ref{i3}), one can write
\begin{equation*}
\widetilde{\phi }_{y^{\prime }}(x^{\prime })=\frac{1}{\sqrt{1-x^{\prime 2}}}%
\frac{1}{\pi }P.V.\int_{-1}^{1}\sqrt{1-\xi ^{2}}\frac{\widetilde{\phi }%
_{x^{\prime }}(\xi )}{x^{\prime }-\xi }d\xi,
\end{equation*}%
which is equivalent, by (\ref{expan}), (\ref{i2}) and (\ref{i4}), to%
\begin{equation*}
\widetilde{\phi }_{y^{\prime }}(x^{\prime })=\partial _{x^{\prime }}\left(
\sqrt{1-x^{\prime 2}}\frac{1}{\pi }P.V.\int_{-1}^{1}\frac{1}{\sqrt{1-\xi ^{2}%
}}\frac{\widetilde{\phi }(\xi )}{x^{\prime }-\xi }d\xi \right) .
\end{equation*}

If we expand, instead of (\ref{expphix}), in the form%
\begin{align}
\widetilde{\phi }_{x^{\prime }} &=\sum_{n=1}^{\infty }a_{n}T_{n}(x^{\prime
}),  \label{ff1} \\
\widetilde{\phi }_{y^{\prime }} &=\sum_{n=1}^{\infty }b_{n}\sqrt{%
1-x^{\prime 2}}U_{n-1}(x^{\prime }),  \label{ff2}
\end{align}%
with $a_{n}$ such that
\begin{equation*}
\int_{-1}^{1}\frac{|\widetilde{\phi }_{x^{\prime }}(x^{\prime
})|^{2}}{\sqrt{1-x^{\prime 2}}}dx^{\prime }=\pi a_{0}^{2}+\frac{%
\pi }{2}\sum_{n=1}^{\infty }a_{n}^{2}<\infty,
\end{equation*}%
then, $b_{n}=a_{n}$ and using (\ref{i2}) one can write
\begin{equation*}
\widetilde{\phi }_{y^{\prime }}(x^{\prime })=\sqrt{1-x^{\prime 2}}\frac{1}{%
\pi }P.V.\int_{-1}^{1}\frac{1}{\sqrt{1-\xi ^{2}}}\frac{\widetilde{\phi }%
_{x^{\prime }}(\xi )}{x^{\prime }-\xi }d\xi.
\end{equation*}

Note that the mass conservation condition $\int_{-1}^{1}\widetilde{\phi }%
_{y^{\prime }}(x^{\prime })dx^{\prime }=0$ implies $b_{1}=0$ which yields $%
a_{1}=0$ implying
\begin{equation*}
\int_{-1}^{1}\frac{x^{\prime }\widetilde{\phi }_{x^{\prime }}(x^{\prime
})dx^{\prime }}{\sqrt{1-x^{\prime 2}}}=0\text{,}
\end{equation*}%
and the absence of $T_{0}(x^{\prime })$ term in (\ref{ff1}) implies%
\begin{equation*}
\int_{-1}^{1}\frac{\widetilde{\phi }_{x^{\prime }}(x^{\prime })dx^{\prime }}{%
\sqrt{1-x^{\prime 2}}}=0\text{.}
\end{equation*}

\section{The compact embedding of $H_{w^{-1}}^{\frac{1}{2}}$ into $%
L_{w^{-1}}^{2}$}

\label{app2}

Since the spaces $H_{w^{-1}}^{\frac{1}{2}}$ and $L_{w^{-1}}^{2}$ can be
defined in terms of an orthogonal basis, following \cite{bisgard2012compact}%
, let us prove the more general result $h^{1/2}\overset{c}{\hookrightarrow}%
\ell^2$ for sequence spaces
\begin{equation*}
h^{1/2}:=\left\{(a_n):\sum_{n=1}^{\infty}n|a_n|^2<\infty\right\}, \quad
\ell^{2}:=\left\{(a_n):\sum_{n=1}^{\infty}|a_n|^2<\infty\right\}.
\end{equation*}

Consider $M\subseteq h^{1/2}$ a bounded set in $h^{1/2}$. Then, given $%
\mathbf{a\in }\mathbb{\ell }^{2}$,
\begin{equation*}
\left\Vert\mathbf{a}\right\Vert_{\mathbb{\ell}^{2}}^{2}=\sum_{n=1}^\infty
a_{n}^{2}\leq \sum_{n=1}^{\infty }na_{n}^{2}=\left\Vert \mathbf{a}
\right\Vert_{h^{\frac{1}{2}}}^{2}<K^{2},
\end{equation*}
where $K$ is such that $\left\Vert \mathbf{a}\right\Vert _{h^{\frac{1}{2}%
}}<K $, $\forall\mathbf{a}\in M$. We take now $Y_{\varepsilon}=span\left\{%
\mathbf{e}_{1},\mathbf{e}_{2},...,\mathbf{e}_{N}\right\}$ where $N$ will be
chosen later as a function of $\varepsilon$, and $\mathbf{e}_{i}$ is the
vector in $\mathbb{\ell}^{2}$ with all components zero except the $i$-th
component that is $1$. We write $\mathbf{x}=a_{1}\mathbf{e}_{1}+a_{2}\mathbf{%
\ e}_{2}+...+a_{N}\mathbf{e}_{N}$ so that $\mathbf{x}\in
Y_{\varepsilon}\subset\mathbb{\ell}^{2}$, $\mathbf{x}%
=(a_{1},a_{2},...,a_{N-1},0,0,...)$ and
\begin{equation*}
\left\Vert \mathbf{x}-\mathbf{a}\right\Vert_{\mathbb{\ell }
^{2}}^{2}=\sum_{n=N}^{\infty}a_{n}^{2}=\sum_{n=N}^{\infty }\frac{1}{n}
na_{n}^{2}\leq\frac{1}{N}\sum_{N}^{\infty}na_{n}^{2}\le\frac{K^{2}}{N}<\frac{
\varepsilon}{2},
\end{equation*}
which implies $N=O(1/\varepsilon )$. Let $\mathbf{x\in }\overline{M}$. There
exists $\mathbf{x}_{n}\in M$ such that $\mathbf{x}_{n}\rightarrow \mathbf{x}$
in the $\mathbb{\ell}^{2}$ topology. For $n$ sufficiently large, $\left\Vert
\mathbf{x}_{n}-\mathbf{x}\right\Vert _{\mathbb{\ell }^{2}}<\frac{\varepsilon
}{2}$ and all $\mathbf{x}_{n}$ are at distance $\frac{\varepsilon}{2}$ of $%
Y_{\varepsilon}$; then by the triangle inequality $\mathbf{x}$ is within $%
\varepsilon$ of $Y_{\varepsilon}$ and the closure of $M$ in $\mathbb{\ell}%
^{2}$ is compact by Proposition 7.4 in \cite{deimling2010} or Proposition
2.1 in \cite{bisgard2012compact}. Hence $h^{1/2}$ is compactly embedded in $%
\mathbb{\ell }^{2}$.

%\section*{Acknowledgments}
%We would like to acknowledge

\bibliographystyle{abbrv}
\bibliography{references}

\end{document}